\documentclass[a4paper,12pt,twoside, fleqn]{amsart}
\usepackage{amsmath,amsthm,amsfonts}
\usepackage[T1]{fontenc}
\newtheorem{thm}{Theorem}[section]

\newtheorem{lem}[thm]{Lemma}

\newtheorem{prop}[thm]{Proposition}
\newtheorem{proposition}[thm]{Proposition}
\newtheorem{definition}[thm]{Definition}

\theoremstyle{definition}

\newtheorem{rem}[thm]{Remark}
\newtheorem{rems}[thm]{Remarks}

\newcommand{\R}{{\mathbb{R}}}

\newcommand{\E}{{\mathbb{E}}}

\newcommand{\C}{{\mathbb{C}}}

\newcommand{\F}{{\mathbb{F}}}
\newcommand{\Z}{{\mathbb{Z}}}
\newcommand{\N}{{\mathbb{N}}}
\newcommand{\T}{{\mathbb{T}}}
\newcommand{\PP}{{\mathbb{P}}}

\newcommand{\Cal}{\mathcal}
\newcommand{\cal}{\mathcal}

\def \Var {{\rm Var}}

\def\j{{\underline j}} \def\k{{\underline k}}

 \def\el{{\underline \ell}}
\def\t{{\underline t}} \def\n {{\underline n}}
\def\0{{\underline 0}}  \def\1{{\underline 1}}

 \def\p {{\underline p}}
\def\u{{\underline u}} 
\def\v{{\underline v}} 
\def\s{{\underline s}} 
 \def\q{{\underline q}}

\def\al {{\underline \alpha}}

\def \and{\text{ and }}

\def \stm0{{\setminus \{\0\}}}

\def \eop{\qed}
\def \Proof {\vskip -3mm {{\it Proof}. }}
\def \proof {\vskip -3mm {{\it Proof}. }}
\def \Proof {\vskip -3mm {{\it Proof}. \ }}
\def \Var  {{\rm Var }}

\def \Card {{\rm Card}}

\begin{document}

\date{\today}

\parskip=2mm
\baselineskip 15pt
\parindent=0mm

\title[On the quenched  functional CLT in 2d random sceneries] {On the quenched functional CLT \\ in 2d random sceneries, examples}
\vskip 2mm

\author{Guy Cohen and Jean-Pierre Conze}
\address{Guy Cohen, \hfill \break Dept. of Electrical Engineering, \hfill \break Ben-Gurion University, Israel} \email{guycohen@bgu.ac.il}
\address{Jean-Pierre Conze, \hfill \break IRMAR, CNRS UMR 6625, \hfill \break University of Rennes I, Campus de Beaulieu, 35042 Rennes Cedex, France} \email{conze@univ-rennes1.fr}

\maketitle

\vskip 3mm 
\begin{abstract}  
We prove a quenched functional central limit theorem (quenched FCLT) for the sums of a random field (r.f.) along a 2d-random walk in different situations: 
when the r.f. is iid with a second order moment (random sceneries), or when it is generated by the action of commuting automorphisms of a torus.
We consider also a quenched version of the FCLT when the random walk is replaced by a Lorentz process in the random scenery.
\end{abstract}

\tableofcontents

\section*{\bf Introduction} 

Let $X = (X_\el)_{\el \in \Z^d}$, $d \geq 1$, be a strictly stationary real random field (r.f.), where the $X_\el$'s have zero mean and finite second moment. 
The r.f. can be represented in terms of dynamical system as $X_\el = T^\el f$, where $(E, \cal A, \mu)$ is a probability space, $T_1, ..., T_d$ are
commuting measure preserving maps on $(E, \cal A, \mu)$ and $f$ is in $L^2(E, \cal A, \mu)$.\footnote{Underlined letters represent elements of $\Z^d$ or $\T^d$. 
We write $\el =(\ell_1, ..., \ell_d)$ and $T^\el f(x) = f(T_1^{\ell_1}... T_d^{\ell_d}x)$. The $\ell^1$-norm in $\Z^d$ of a vector $\el$ is denoted $|\el|$.}
\footnote{If the maps are not invertible, using stationarity the random field can be extended to a strictly stationary random field indexed by $\Z^d$.} 

Let $w= (w_n)_{n \geq 1}$ be a sequence of weights (or summation sequence), that is for each $n$ a function $\el \in \Z^d \to w_n(\el) \in \R$, 
with $0 < \sum_{\el \in \Z^d} |w_n(\el)| < +\infty$. A natural question is the asymptotic normality in distribution of the self-normalized sums 
$\displaystyle \sum_{\el \in \Z^d} w_n(\el) \, f(T^\el x) / \|\sum_{\el \in \Z^d} w_n(\el) \, T^\el f\|_2$ 
and the estimation of the normalization factor. A stronger property, for some models, is the validity of a functional central limit theorem (FCLT).

Previously (\cite{{CohCo17}, CohCo16}), we have considered quenched central limit theorems for summation along a random walk, 
as well as summation on a sequence of sets in $\Z^d$. In a forthcoming paper, the FCLT for summation over sets will be presented. 
 
The present paper is about the random walk case and specially the 2-dimensional random walk, the case of $d$-dimensional random walks being easier for $d > 2$.
We show a FCLT in different models for the sums along a r.w. for almost all realizations of the r.w. (quenched FCLT).

One of these models is the random walk in random sceneries, i.e., the sums along a r.w. of a 2-d random field of iid r.v.s with a moment of order 2. 
This improves a result of \cite{GuPoRS14} which uses a slightly stronger moment condition. Our proof is short and self-contained. 
The same method can be used when the usual random walk is replaced by a plane Lorentz process (generated by a periodic billiard with dispersive obstacles) as in \cite{Pen09}. 
A key step in the proof is then the law of large numbers shown in \cite{Pen14} for the self-intersection of the billiard map.
The random sceneries can be also replaced by a random field which is no more iid, but generated by an algebraically defined $\Z^2$-dynamical system.
In this framework, we consider algebraic actions on tori by commuting automorphisms.

Tightness of the process is one of the main step of the proof of a FCLT. In the framework of sums along a random walk, our purpose is to present two different situations, 
independent case and algebraic case, as an illustration of two methods: one relying on the maximal inequality for associated r.v.s as shown by Newman and Wright \cite{NeWr81}, 
the other on norm estimates for the maximum of partial sums as in Billingsley \cite{Bi99}, M\'oricz \cite{Mo76} and others authors.

In Section \ref{varGene1} we gather results about the variance for the sums along a random walk. The independent case is presented in Section \ref{indepSec}.
Some facts on cumulants are recalled in Section \ref{Sectcumulant}, then applied to moving averages in Section \ref{movAverag} and to algebraic models 
in Section \ref{algMod1}.  For the tightness in the latter case, we use a method based on a moment inequality for maxima, which is presented in Section \ref{SectMoricz}.

\section{\bf Summation along a r.w. and variance} \label{varGene1}

Everywhere we assume (or prove) the absolute summability of the series of decorrelations:
\begin{eqnarray}
\sum_{\k \in \Z^d} | \int_X T^\k f \, \overline f \, d\mu| < \infty, \label{sumcorr1}
\end{eqnarray}
an hypothesis which implies existence and continuity of the spectral density associated to $f$, i.e., existence of a function $\varphi_f \in C(\T^d)$ such that
\begin{eqnarray}
\int_X T^\k f \, \overline f \, d\mu = \int_{\T^d} e^{2\pi i \langle\k, \t\rangle} \, \varphi_f \, d\t, \ \forall \k \in \Z^d. \label{specdens1}
\end{eqnarray}

A method of summation is given by random walks (r.w.). If $(Z_n)_{n \geq 0}$ is a random walk starting from $\0$ on $\Z^d$,
the associated ``ergodic sums'' along the orbits of the random walk are  
\begin{eqnarray}
\sum_{k=0}^{n-1} \, T^{Z_k(\omega)} f =  \sum_{\el \in \Z^d} w_n(\omega, \el) \, T^\el f, \text{ with } w_n(\omega, \el) = \#\{k < n: Z_k(\omega) = \el\}.
\end{eqnarray}

{\it Remark.} Summation along the orbits of the random walk differs from summation over the range of the random walk.  
It has been shown (\cite{DeKo17}) that the range of the random walk $(Z_n)$ has the F\o{}lner  property, so that summation over the range of the random walk 
yields a CLT. Nevertheless, a functional CLT for summation over the range is a question. 

\vskip 2mm
\subsection{\bf Variance} \label{varSect}

First let us recall some results on the variance $\int_E |\sum_{\el \in \Z^d} w_n(\el) \, T^\el f|^2 \, d\mu$ which will be useful for the FCLT. 
Its computation is based on the normalized non-negative kernel
\begin{eqnarray}
K(w_n)(\t) = {|\sum_{\el \in \Z^d} w_n(\el) \, e^{2\pi i \langle \el, \t \rangle}|^2 \over \sum_{\el \in \Z^d} |w_n(\el)|^2}, \ \t \in \T^d.\label{chch}
\end{eqnarray}
We say that $w = (w_n)$ is $\xi$-{\it regular}, where $\xi$ is a probability measure on $\T^d$, if $(K(w_n))_{n \geq 1}$ weakly converges to $\xi$: 
$\lim_{n \to \infty}\int_{\T^d} K(w_n) \,\varphi \, d\t  = \xi(\varphi)$ for every continuous $\varphi$ on $\T^d$,
or equivalently if 
\begin{eqnarray*}
\hat \xi (\p)  = \lim_{n \to \infty} \int K(w_n)(\t) \, e^{-2\pi i \langle \p, \t \rangle} \ d\t, \, \forall \p \in \Z^d.
\end{eqnarray*}
Under Condition (\ref{sumcorr1}), the asymptotic variance for $f$ is then
\begin{eqnarray}
\sigma_w^2(f) := \lim_n \bigl(\|\sum_{\el \in \Z^d} w_n(\el) \, T^\el f\|_2^2/ \sum_{\el \in \Z^d} |w_n(\el)|^2 \bigr)  = \xi(\varphi_f). \label{asymptVar1}
\end{eqnarray}
In what follows, we will deal with examples which are $\delta_\0$-regular. For examples of summation along a random walk which are $\xi$-regular 
with $\xi \not = \delta_0$, see for instance \cite{CohCo17}.

\begin{rem} \label{densRem1} If $\xi$ is a probability measure on $\T^d$, $f \to (\xi(\varphi_f))^\frac12$ satisfies the triangular inequality.
Indeed, for $p(\t) = \sum a_\el \, e^{2\pi i \langle \el, \t\rangle}$ a trigonometric polynomial, we have: 
$\int |p(\t)|^2 \, \varphi_f(\t) \, d\t = \|\sum a_\el \, T^\el f \|_2^2$, by definition of the spectral density; hence, by the triangular inequality,
$$(\int |p(\t)|^2 \, \varphi_{f+g}(\t) \, d\t)^\frac12 \leq (\int |p(\t)|^2 \, \varphi_{f}(\t) \, d\t)^\frac12
+ (\int |p(\t)|^2 \, \varphi_{g}(\t) \, d\t)^\frac12.$$
It follows $\varphi_{f+g}^\frac12 \leq \varphi_{f}^\frac12  + \varphi_{g}^\frac12$, if $\varphi_{f+g}, \varphi_{f}, \varphi_{g}$ are continuous;
hence: $\varphi_{f+g} \leq \varphi_{f}  + \varphi_{g} + 2 \varphi_{f}^\frac12 \varphi_{g}^\frac12$, which implies
$\xi(\varphi_{f+g}) \leq \xi(\varphi_{f})  + \xi(\varphi_{g}) + 2 \xi(\varphi_{f}^\frac12 \varphi_{g}^\frac12)
\leq \xi(\varphi_{f})  + \xi(\varphi_{g}) + 2 (\xi(\varphi_{f}))^\frac12 \, (\xi(\varphi_{g}))^\frac12$.
\end{rem}

{\it Question of non-degeneracy}

Condition (\ref{sumcorr1}) implies that, for any $\delta_\0$-regular summation sequence, the following conditions are equivalent:  nullity of the asymptotic variance,
$\varphi_f(\0) = 0$, $\sum_{\k \in \Z^d} \langle T^\k f \, f \rangle = 0$.

A function $f$ is called a {\it mixed coboundary} if there exists measurable functions $g_i$, $i=1, ..., d$, such that $f = \sum_{i=1}^d (I - T_i) g_i$. 
In the example of a $\Z^d$-action by commuting algebraic automorphisms of a torus, for a class of regular function, the nullity of the asymptotic variance occurs 
if and only if $f$ is a mixed coboundary. (See \cite{CohCo13}). 

\vskip 3mm
\subsection{\bf Sums along random walks} \label{varRWSect1}

\

{\it Definitions and notations.}

Let $(\zeta_i, \, i = 0, 1, ...)$ be a sequence of i.i.d. random vectors on a probability space $(\Omega, \, \PP)$ with values in $\Z^d$.
The corresponding {\it random walk} (r.w.) $Z = (Z_n)$ in $\Z^d$ starting from $\0$ is defined by $Z_0 := \0$, $Z_n := \zeta_0 +... + \zeta_{n-1}$, $n \geq 1$.
We suppose $Z$ to be aperiodic
\footnote{i.e., we suppose that the subgroup generated in $\Z^d$ by $\{\el: \PP(\zeta_\0 = \el) > 0\}$ is $\Z^d$.}, 
with 0 mean, finite variance and (nonsingular) covariance matrix $\Sigma$. (For random walks, see \cite{Sp64}.)

The r.v.s $\zeta_i$ can be viewed as the coordinate maps on $(\Omega, \, \PP)$ obtained as  $(\Z^d)^\Z$ equipped with a product measure and with the shift $\theta$ acting on the coordinates. 
We have $\zeta_i = \zeta_0 \circ \theta^i$ and the cocycle relation $Z_{n +n'} = Z_n + Z_{n'} \circ \theta^n, \forall n, n' \in \Z,$ holds.

Given a random field\footnote{Recall that the process is denoted either by $(X_\el)$ or by $(T^\el f)$.}  $X = (X_\el$, $\el \in \Z^d)$ on $(E, \mu)$, we form the process 
on $(E, \mu)$ obtained by summing along the r.w. $(Z_n)$. We denote these sums, for a fixed $\omega$, by 
\begin{eqnarray}
&&S_n^{\omega, X}(x) = S_n^{\omega}(x) :=\sum_{i=0}^{n-1} X_{Z_i(\omega)}(x), \, n \geq 1. \label{defsumRW}
\end{eqnarray}
If the random field is represented as $X_\el = T^\el f$, the sums read:
\begin{eqnarray}
S_n^\omega f = \sum_{i=0}^{n-1} T^{Z_i(\omega)} \, f = \sum_{\el \in \Z^d} w_n(\omega, \el) \, T^{\el} \, f, 
\text{ with } w_n(\omega, \el) = \sum_{0 \leq i < n} 1_{Z_i(\omega) = \el}. \label{defwn0}
\end{eqnarray}
Summing along the random walk amounts to fix $\omega$ in the ergodic sums of the ``skew product'':
$(\omega, x) \to T_{\zeta_0}(\omega, x) = (\theta \omega, T^{\zeta_0(\omega)} x)$ on $\Omega \times E$. Putting 
$F(\omega, x) = f(x)$, for an observable $f$ on $E$, we get that the ergodic sums of $F$ for $T_{\zeta_0}$ read:
\begin{eqnarray}
S_n F (\omega, x) = \sum_{i=0}^{n-1} F(T_{\zeta_0}^i \, (\omega, x)) =  \sum_{i=0}^{n-1} f(T^{Z_i(\omega)} x)
= (S_n^\omega f)(x). \label{ergSum1}
\end{eqnarray}
If we consider the r.v. $S_n F (\omega, x)$ as defined on $\Omega \times E$ endowed with the probability $\PP \times \mu$, a limit theorem is sometimes 
called {\it annealed}. We can also fix $\omega \in \Omega$. 
A limit theorem in distribution (with respect to the measure $\mu$ on $E$) obtained for $\PP$-a.e. $\omega$ is called {\it quenched}.

{\it We will consider the case where $Z$ is a r.w. in $\Z^2$}. In this case, $(Z_n )$ is recurrent and a non standard normalization occurs in the CLT for sums along $Z_n$ as recalled below.

\vskip 3mm
\subsubsection{\bf On the number of self-intersections of a r.w.}

\

If $I, J$ are intervals, the quantity $V(\omega, I, J, \p) =$
\begin{eqnarray}
&&\int \bigl(\sum_{u \in I}  e^{2\pi i \langle Z_u(\omega), \t\rangle}\bigr) 
\, \bigl(\sum_{v \in J} \, e^{- 2\pi i \langle Z_v(\omega), \t \rangle}\bigr) \, e^{-2\pi i \langle \p, \t\rangle} \, d\t 
\ = \#\{(u, v) \in I \times J: \, Z_u(\omega) - Z_v(\omega) = \p\} \label{cross1}
\end{eqnarray}
is non negative and increases when $I$ or $J$ increases for the inclusion order.

We write simply $V(\omega, I, \p)$ for $I=J$,  $V(\omega, I)$ for  $V(\omega, I, \0)$ and $ V_n(\omega)$ for $V(\omega, [0, n[)$.
Hence $V_n(\omega) = \#\{0 \leq u, v < n: \ Z_u(\omega) = Z_v(\omega)\}$ is the number of self-intersections starting from $\0$.

Observe that $V(\omega, J) = \sum_{\el \in \Z^2} w(\omega, J, \el)^2$, where $w_n(\omega, J, \el) = \sum_{i \in J} 1_{Z_i(\omega) = \el}$. In particular $V_n(\omega) = \sum_{\el \in \Z^2} w_n(\omega, \el)^2$.

Note also that $V(\omega, [b, b+k[) = V(\theta^b\omega, [0, k[) = V_k(\theta^b\omega)$, for $b \geq 0, k \geq 1$.

Let $A, B$  be  in $[0, 1]$. We have \footnote{For simplicity, in the formulas below, we write $nA$, $nB$ instead of $\lfloor nA \rfloor$ or $\lfloor nA \rfloor +1$, 
$\lfloor nB \rfloor$, $\theta^t$ instead of $\theta^{\lfloor t \rfloor}$. The equalities are satisfied up to the addition of quantities 
which are bounded independently from $A, B, n$.}
\begin{flalign*}
&V(\omega, [nA, nB], \p) =  \int \bigl[(\sum_{u \in [nA, nB]}  e^{2\pi i \langle Z_u(\omega), \t\rangle}) 
\, (\sum_{v\in [nA, nB]} \, e^{- 2\pi i \langle Z_v(\omega), \t \rangle}) \, e^{-2\pi i \langle \p, \t\rangle} \bigr] \, d\t& \\
&= \#\{u, v \in [0, n(B-A)]:\, \sum_{i=0}^{u-1} \zeta(\theta^{i+nA} \omega) - \sum_{i=0}^{v-1} \zeta(\theta^{i+nA} \omega) = \p\}
= V(\theta^{nA}\omega, [0, n(B-A)], \p).&
\end{flalign*}
For $d=2$, there are $C_0, C$ finite positive constants
\footnote{If the r.w. is strongly aperiodic, $C_0 = (\pi \sqrt{\det \Sigma})^{-1}$. For a general aperiodic r.w. in $\Z^2$, see for instance Theorem 5.1 in \cite{CohCo17}.}
such that the following laws of large numbers hold 
(see: \cite{Bo89} Lemma 2.6 for (\ref{BoundVar}), \cite{Lew93} step 1 in the proof of Proposition 1.4 for (\ref{equivAn20}), 
and \cite{CohCo17} Theorem 3.13 for (\ref{equivAn2})):
\begin{eqnarray}
&&\E(V_n) \sim C_0 n \ln n, \ \Var(V_n) \leq C n^2, \label{BoundVar}\\
&&\varphi_n(\omega) := {V_n(\omega) \over C_0 n \ln n} \to 1, \text{ for a.e. } \omega, \label{equivAn20} \\ 
&&\varphi_n(\omega, \p) := {V(\omega, [1,n], \p) \over C_0 n \ln n}  \to 1, \forall \p \in\Z^d, \text{ for a.e. } \omega. \label{equivAn2}
\end{eqnarray}

We denote by $\Omega_0$ the set of full probablity of $\omega$'s for which (\ref{equivAn20}) and (\ref{equivAn2})  hold.
We have
\begin{eqnarray}
&&V_n(\omega) \leq K(\omega) \, n \ln n, \, \forall n \geq 2, \text{ where the function } K \geq 0 \text{ is finite on } \Omega_0, \label{Komega}\\
&& \text{ for any fixed } A \in ]0, 1], \ V(\omega, [1, nA], \p) \sim C_0 n A \ln n, \text{ for } \omega \in \Omega_0. \label{equivAn}
\end{eqnarray}

Recall that (\ref{equivAn2}) shows the $\delta_\0$-regularity of the summation sequence along the random walk $Z$ for a.e. $\omega$ (cf. \cite{CohCo17}): 
if $f$ has a continuous spectral density $\varphi_f$,
\begin{eqnarray}
&&(C_0 n \ln n)^{-1} \|\sum_{k=0}^{n} T^{Z_k(\omega)} f\|_2^2 = (C_0 n \ln n)^{-1} 
\int |\sum_{k=0}^{n} e^{2\pi i \langle Z_k(\omega) , \t \rangle}|^2 \, \varphi_f(\t) \, d\t \to \varphi_f(\0). \label{omegaReg1}
\end{eqnarray}

Before a preliminary lemma, let us introduce some more notations. For $J \subset \N$, we put:
\begin{eqnarray*}
&&U^{(m)}(\omega, J) := \sum_\el w(\omega, J, \el)^m, \ \ U_n^{(m)}(\omega) := \sum_\el w_n(\omega, \el)^m.
\end{eqnarray*}
For $\el_1, \el_2, \el_3 \in \Z^d$, we put $W_n(\omega, \el_1, \el_2, \el_3) :=$
\begin{equation} 
\#\{1 \leq i_0, \, i_1, \, i_2, \, i_3 < n: \, Z_{i_1}(\omega) - Z_{i_0}(\omega) = \el_1, 
Z_{i_2}(\omega) - Z_{i_0}(\omega)= \el_2, Z_{i_3}(\omega) - Z_{i_0}(\omega) = \el_3\}. \label{wn123}
\end{equation}
By \cite[Lemma 2.5]{Bo89} (see also \cite[Proposition 2.9] {CohCo17}) we have, for every $\varepsilon>0$, 
\begin{eqnarray}
\sup_{\el\in \Z^2}w_n(\omega, \el)=o(n^\varepsilon), \text{ for a.e. } \omega. \label{nVarep1}
\end{eqnarray}
Therefore, for every $\varepsilon>0$, there is $C(\omega) = C(\omega, m, \varepsilon)$ such that
\begin{eqnarray}
U_n^{(m)}(\omega)  = \sum_\el w_n(\omega, \el)^m \leq C(\omega) \, n^{1+\varepsilon} \label{majwnmIneg0}.
\end{eqnarray}
\begin{lem} \label{majwnm} 
There exists a positive integrable function $C_3$ such that 
\begin{eqnarray} 
W_n(\omega, \el_1, \el_2, \el_3) \leq C_3(\omega) \, n \, (\ln n)^5, \, \forall n \geq 1. \label{majnumb3}
\end{eqnarray}
\end{lem}
\proof Since the terms in the sum (\ref{wn123}) with equality between indices  can be treated by induction, it suffices to bound 
$$W_n'(\omega) = \sum_{1 \leq i_0 < i_1 < i_2 < i_3 \leq n} 1_{Z_{i_1} - Z_{i_0} = \el_1} . 1_{Z_{i_2} - Z_{i_1}= \el_2 - \el_1}. 
1_{Z_{i_3} - Z_{i_2} = \el_3 - \el_2}.$$
Using independence and the local limit theorem for the random walk, we find the bound
\begin{eqnarray} 
\int W_n'(\omega) \, d\PP(\omega)  \leq  C_1 \sum_{i_0, i_1, i_2, i_3 \in [1, n]} (i_1 \, i_2 \, i_3)^{-1} \leq C_2 \, n \, (\ln n)^3. \label{majnlgn3}
\end{eqnarray}
By (\ref{majnlgn3}) and similar bounds for the others configurations, we have $\int W(n) \, d\PP \leq C_2' n \, (\ln \, n)^3$ for another constant $C_2'$. 
Therefore $\displaystyle \sum_{p=1}^\infty \int 2^{-p}(\ln(2^p))^{-5} \, W_{2^p} \, d\PP < \infty$.
The function $C(\omega) := \sum_{p=1}^\infty 2^{-p}(\ln(2^p))^{-5} \, W_{2^p}$ is integrable and we have:
$W_{2^p}(\omega) \leq C(\omega) \, 2^p \, (\ln(2^p))^5, \, \forall p \geq 1$.

Let $p_n$ be such that: $2^{p_n-1} \leq n < 2^{p_n}$. Since $W_n$ is increasing with $n$, we obtain:
$$W_{n}(\omega) \leq W_{2^{p_n}}(\omega) \leq K(\omega) \, 2^{p_n} \, (\ln 2^{p_n})^{5} \leq K(\omega) \, 2n (\ln 2n)^{5}
\leq K'(\omega) \, n (\ln n)^{5}. \eop$$ 
The same method shows that, for $m \geq 1$, there is a positive integrable function $K_m$ such that
\begin{eqnarray}
U_n^{(m)}(\omega)  = \sum_\el w_n(\omega, \el)^m \leq K_m(\omega) \, n \, (\ln n)^{m+1}, \, \forall n \geq 1. \label{majwnmIneg1}
\end{eqnarray}

\vskip 2mm
{\bf Study of the variance for the finite dimensional distributions}

The following lemma will be applied to the successive return times of a point $\omega$ into a set under the iteration of the shift $\theta$.
\begin{lem} \label{mean0} Let $(y(j), j \geq 1)$ be a sequence with values in $\{0, 1\}$ such that $\lim_n \frac1n \sum_{j=1}^n y(j) = a > 0$.
If $(k_r)$ is the sequence of successive times such that $y(k_r) = 1$, then, for every $\delta >0$, there is $n(\delta)$ such that,
 for $n \geq n(\delta)$, $k_{r+1} - k_r \leq \delta n$, for all $r \in [1, n]$.
\end{lem}
\Proof Since $r = \sum_{j=1}^{k_r} \, y(j)$, we have: $k_r / r = k_r / \sum_{j=1}^{k_r} \, y(j)  \to a^{-1}$. 
Hence, for every $\delta > 0$, there is $n_1(\delta)$ such that $0 < k_{r+1}  - k_r \leq \delta r$, for $r \geq n_1(\delta)$.
Therefore, if $n \geq n_1(\delta)$, then  $0 < k_{r+1}  - k_r \leq \delta r \leq \delta n$, for $r \in [n_1(\delta), n]$. 

If $n(\delta) \geq n_1(\delta)$ is such that $k_{r+1}  - k_r \leq \delta n(\delta)$ for $r \leq n_1(\delta)$, we get the result of the lemma. \eop

\begin{lem} \label{mean1} Let $\Lambda$ be a measurable set in $\Omega$ of positive measure. 
Let $k_r= k_r(\omega)$ be the successive times such that $\theta^{k_r} \omega \in \Lambda$.
For a.e. $\omega$, for every positive small enough $\delta$, there is $n(\delta)$ such that for $n \geq n(\delta)$
\hfill \break
1) $k_{r+1} - k_r \leq \delta n$, for all $r \in [1, n]$; moreover, $k_n \sim c n$, where $c=  \PP(\Lambda)^{-1}$;
\hfill \break
2) there are integers $v < 2/\delta$ and $0 = \rho_1^{(n)} < \rho_2^{(n)} < ... < \rho_v^{(n)} \leq n < \rho_{v+1}^{(n)}$, 
such that $\theta^{\rho_i^{(n)}} \omega \in \Lambda$ and $\frac12 \delta n \leq \rho_{i+1}^{(n)} - \rho_i^{(n)} \leq \frac32 \delta n$, for $i = 1, ..., v$.
\end{lem}
\Proof Since $\theta$ is ergodic on $(\Omega, \PP)$, Birkhoff ergodic theorem implies
$\lim_n \frac1n \sum_0^{n-1} 1_\Lambda (\theta^k \omega) = \PP(\Lambda) > 0$, for a.e. $\omega$ and $k_n / n \to \PP(\Lambda)^{-1}$.
Hence Lemma \ref{mean0} implies 1). For 2), we select an increasing sequence of visit times to the set $\Lambda$ satisfying the prescribed conditions by eliminating successive 
times at a distance $< \frac12 \delta n$. \eop

\vskip 2mm
{\it Asymptotic orthogonality of the cross terms}
 
We show the asymptotic orthogonality of the cross terms:  for $0 < A < B < C < D < 1$, $\p \in \Z$,
\begin{eqnarray}
&&\int \bigl[(\sum_{v = nA}^{nB}  e^{2\pi i \langle Z_v(\omega), \u \rangle}) 
\, (\sum_{w = nC}^{nD} \, e^{- 2\pi i \langle Z_w(\omega), \u \rangle}) \, e^{-2\pi i \langle \p, \u \rangle} \bigr] \, d\u
=\varepsilon_n(\omega) \, n \log n, \text{ with } \varepsilon_n(\omega) \to 0. \label{orth1}
\end{eqnarray}
The above integral is the non negative self-intersection quantity: $V(\omega, [nA, nB], [nC, nD], \p)$. By (\ref{cross1}), $V(\omega, I, J, \p)$ increases when $I$ or $J$ increases. 
Hence, it suffices to show (\ref{orth1}) for the intervals $[1, nA], [nA, n]$, for $0 < A < 1$. The proof below is based on (\ref{equivAn2}) and (\ref{equivAn}).

\begin{lem} \label{asympOrtho1} There is a set $\hat \Omega \subset \Omega$ such that $\PP(\hat \Omega) = 1$ and for all $\omega \in \hat \Omega$,
the following holds:
\begin{eqnarray}
&&\lim_n \varphi_{n B}(\theta^{n A}\omega, \, \p) = \lim_n {V(\omega, [nA, n], \p) \over C_0 \, n B \ln n} = 1, \text{ for } A\in ]0,1[, B = 1 - A; \label{lowup0} \\
&&V(\omega, [1, nA], [nA, n], \p) + V(\omega,[nA, n], [1, nA], \p) =\varepsilon_n(\omega) \, n \log n, \text{ with } \varepsilon_n(\omega) \to 0. \label{orth2}
\end{eqnarray}
\end{lem}
\proof  {\it 1) The set $\hat \Omega$.} 
For every $L \geq 1$ and $\delta > 0$, let $\Lambda(L, \delta):= \{\omega: \varphi_n(\omega, \, \p)  - 1 \in [-\delta, \delta], \forall n \geq L\}$.
We have $\lim_{L \uparrow \infty} \PP(\Lambda(L, \delta)) = 1$. There is $L(\delta)$ such that $\PP(\Lambda(L(\delta), \delta)) \geq \frac12$.

We will apply Lemma \ref{mean1} to $\Lambda(L(\delta_j), \delta_j)$ for each $j$, where $(\delta_j)$ is a sequence tending to 0, 
therefore getting a set $\omega$'s of full $\PP$-measure. The set $\hat \Omega$ is the intersection of this set with the set $\Omega_0$ (of full measure) 
for which the law of large numbers holds for $(V_n(\omega))$. Let $\omega \in \hat \Omega$. 

{\it 2) Proof of (\ref{lowup0})}. 
We have $V(\omega, [nA, n[, \p) = V(\theta^{nA} \omega, [0, n(1-A)[, \p)$ and
\begin{eqnarray}
&&V(\omega, [1, n], \p) - V(\omega, [1, nA[, \p) - V(\omega, [nA, n], \p) \nonumber \\
&& = V(\omega, [1, nA[, [nA, n[, \p) + V(\omega,[nA, n], [1, nA[, \p) \geq 0. \label{differ1}
\end{eqnarray}
Claim: for an absolute constant $C_1$ depending on $A$ and $\p$, for every $\delta$, for $n$ big enough,
\begin{eqnarray}
\varphi_{n B}(\theta^{n A}\omega, \, \p) = {V(\omega, [nA, n], \p) \over C_0 \, n (1 - A) \ln n} \in [1 - C_1 \delta, 1 + C_1 \delta]. \label{lowup1}
\end{eqnarray}

{\it Upper bound:} The law of large numbers for $V_n(\omega, \p)$ implies, with $|\varepsilon_n|, |\varepsilon_n'| \leq \delta$ for $n$ big enough,
$$C_0^{-1} \, V(\omega, [1, n], \p) =  (1 + \varepsilon_n) \, n \ln n, \ C_0^{-1} \, V(\omega, [1, nA], \p) = (1 + \varepsilon_n') \, nA \ln n.$$
With $B= 1 - A$, this implies by (\ref{differ1}) 
\begin{eqnarray*}
{V(\omega,  [nA, n], \p) \over C_0 \, nB \ln n} \leq {(1+ \varepsilon_n) \, n \ln n - (1+ \varepsilon_n') \, n A \ln n \over nB \ln n} 
\leq 1 + {|\varepsilon_n| \over B} + {|\varepsilon_n'| A \over B} \leq 1+ {1 + A \over B} \delta.
\end{eqnarray*}

{\it Lower bound:} We apply Lemma \ref{mean1} to $\Lambda(L(\delta), \delta)$.
 Let $n_A, n_A'$ be two consecutive visit times $\leq n$ such that $n_A \leq n A < n_A'$. 
For $n$ big enough, we have $0 <  n_A' -  n_A \leq \delta n$ and 
$$n_A = n A \, (1 - \rho_n), \ n_A' = n A \, (1 + \rho_n'), \text { with } 0 \leq A \rho_n, A \rho_n'  \leq \delta.$$
Moreover, since $\omega \in \Omega_0$, there is $L$ such that $\varphi_n(\omega, \, \p) - 1 \in [-\delta, +\delta]$ for $n \geq L$.  

We have, with $|\delta_n'| \leq \delta$,
\begin{eqnarray*}
&C_0^{-1} \, V(\omega, [n'_{A}, n], \p) \geq (1-\delta_n') (n - n_A') \,\ln(n- n_A') = (1-\delta_n')  (nB -nA \rho_n') \,\ln(nB -nA \rho_n') .
\end{eqnarray*}
It follows, for $\delta$ (hence $\rho_n'$) small:
\begin{flalign*}
&{V(\omega, [n_A', n], \p) \over C_0\,  (1-\delta_n')  \, n B \ln (n B)} \geq {(nB -nA \rho_n') \,\ln(nB -nA \rho_n') \over n B \ln (n B)} 
= {(B -A \rho_n') \, [\ln(nB) + \ln (1 - \frac AB \rho_n')] \over B \ln (n B)}& \\
&\geq (1- \frac AB  \rho_n') - 2 (1- \frac AB \rho_n')\frac{\frac AB \rho_n'}{\ln(nB)} 
\geq 1- \frac AB \rho_n'- 2 \frac{\frac AB \rho_n'}{\ln(nB)} \geq 1-B^{-1}\delta (1 + \frac{2}{\ln(nB)}).&
\end{flalign*}
As $V(\omega, J, \, \p)$ increases when the set $J$ increases, we have by the choice of $n_A$ and $n_A'$:
$$V(\omega, [n_A', n], \p) \leq V(\omega, [nA, n], \p).$$ 
Therefore, for $n$ such that $\ln (nB) \geq 2$, we have
$${V(\omega, [nA, n], \p) \over C_0 \, n B \ln (n B)} \geq (1-\delta) \, (1- {2 \over B}\delta) \geq 1 - \delta (1 + {2 \over B}).$$
This shows the lower bound. Altogether with the upper bound, this proves the claim (\ref{lowup1}).

{\it 3)  Proof of (\ref{orth2}).} Let $\delta > 0$. According to (\ref{differ1}) and (\ref{lowup1}), for $n$ big enough, we have
\begin{flalign*}
&V(\omega, [1, nA], [nA, n], \p) + V(\omega,[nA, n], [1, nA], \p) = V(\omega, [1, n], \p) - V(\omega, [1, nA], \p) - V(\omega, [nA, n], \p) &\\
&=  C_0[ (1 + \varepsilon_n) \, n \ln n  - (1 + \varepsilon_n') \, n A \ln n - (1 + \varepsilon_n'') \, n(1-A) \ln n  \leq (2 + C_1) \, C_0 \, \delta \, n \ln n.  \eop&
\end{flalign*}

For the asymptotic variance for $\sum_{j= 0}^s \, a_j \sum_{k=nt_{j-1} }^{nt_j} T^{Z_k(\omega)} f$, 
where $a_1, ..., a_s$ are real numbers and $0 = t_0 < t_1 < ... < t_{s-1} < t_s = 1$ is a subdivision of $[0, 1]$, 
we will use the following lemma, to which we will refer for the processes $(S_n^{\omega, X})$ considered later.
\begin{lem} \label{VarrwLem1} Assume that $f$ has a continuous spectral density $\varphi_f$. For a.e. $\omega$ we have
\begin{eqnarray}
(C_0 n \ln n)^{-1} \|\sum_{j= 1}^s \, a_j \sum_{k=nt_{j-1}}^{nt_j} T^{Z_k(\omega)} f\|_2^2
\to \varphi_f(\underline 0)\sum_{j=1}^s a_j^2(t_j-t_{j-1}). \label{VarMulti1}
\end{eqnarray}
\end{lem}
\proof 1) Recall that proving (\ref{VarMulti1}) amounts to prove
$$(C_0 n \ln n)^{-1}  \int |\sum_{j= 1}^s \, a_j \sum_{k=nt_{j-1}}^{nt_j} e^{2\pi i \langle Z_k(\omega) , \u \rangle}|^2 \, \varphi_f(\u) \, d\u
\to \varphi_f(\underline 0)\sum_{j=1}^s a_j^2(t_j-t_{j-1}).$$

1) First suppose that $\varphi_f$ is a trigonometric polynomial $\rho$, which allows to use (\ref{orth1}) for a finite set of characters $e^{-2\pi i \langle \p, \u \rangle}$.
Using (\ref{omegaReg1}) for the asymptotic variance starting from 0, we have 
$(C_0 n \ln n)^{-1} \|\sum_{k=0}^{\lfloor tn \rfloor} T^{Z_k(\omega)} f\|_2^2 \to  t \rho(\0)$, for $t \in ]0, 1[$.
By Lemma \ref{asympOrtho1}, 
\begin{eqnarray*}
&&(C_0 n \ln n)^{-1} \|\sum_{k=\lfloor s n \rfloor}^{\lfloor t n \rfloor} T^{Z_k(\omega)} f\|_2^2 \to  (t-s) \, \rho(\0), \text{ for } 0 < s < t < 1.
\end{eqnarray*}
Expanding the square and using that the cross terms are asymptotically negligible, we have
\begin{eqnarray*}
&&(C_0 n \ln n)^{-1} \int |\sum_{j= 1}^s \, a_j \sum_{k=nt_{j-1}}^{nt_j} e^{2\pi i \langle Z_k(\omega) , \u \rangle}|^2 \, \rho(\u) \, d\u \\
&& \sim \ (C_0 n \ln n)^{-1} \bigl(\sum_{j=1}^s a_j^2 \,  \int |\sum_{k=nt_{j-1}}^{nt_j} e^{2\pi i \langle Z_k(\omega) , \u \rangle}|^2 \, \rho(\u) \, d\u \bigr)
\to \rho(\0) \sum_{j=1}^s a_j^2(t_j-t_{j-1}).
\end{eqnarray*}
This shows (\ref{VarMulti1}) for trigonometric polynomials. 

2) For a general continuous spectral density $\varphi_f$, for $\varepsilon >0$, 
let $\rho$ be a trigonometric polynomial, such that  $\|\varphi_f - \rho\|_\infty < \varepsilon$. Remark that
$$\int |\sum_{j= 1}^s \, a_j \sum_{k=nt_{j-1}}^{nt_j} e^{2\pi i \langle Z_k(\omega) , \u \rangle}|^2 \, d\u 
\leq \sum_{j, j'= 1}^s \,a_j a_{j'} \, V(\omega, [nt_{j-1}, nt_j],  [nt_{j'-1}, nt_{j'}], \0) \leq (\sum_{j= 1}^s \,|a_j|)^2 \, V_n(\omega).$$
Therefore we have:
\begin{eqnarray*}
&&\bigl|(C_0 n \ln n)^{-1} \int |\sum_{j= 1}^s \, a_j \sum_{k=nt_{j-1}}^{nt_j} e^{2\pi i \langle Z_k(\omega) , \u \rangle}|^2 \, \varphi_f(\u) \, d\u
- \varphi_f(\0) \sum_{j=1}^s a_j^2(t_j-t_{j-1})\bigr|\\
&&\leq \bigl|(C_0 n \ln n)^{-1} \int |\sum_{j= 1}^s \, a_j \sum_{k=nt_{j-1}}^{nt_j} e^{2\pi i \langle Z_k(\omega) , \u \rangle}|^2 \, \rho(\u) \, d\u
- \rho(\0) \sum_{j=1}^s a_j^2(t_j-t_{j-1})\bigr| \\
&&+ \varepsilon \, [(C_0 n \ln n)^{-1} \int |\sum_{j= 1}^s \, a_j \sum_{k=nt_{j-1}}^{nt_j} e^{2\pi i \langle Z_k(\omega) , \u \rangle}|^2 \, d\u
+ \sum_{j=1}^s a_j^2(t_j-t_{j-1})].
\end{eqnarray*}
By the remark, the above quantity inside $[\ ]$ is less than $(\sum_{j= 1}^s \,|a_j|)^2 \, C_0 n \ln n)^{-1} V_n(\omega) + \sum_{j=1}^s a_j^2(t_j-t_{j-1})$,
which is bounded uniformly with respect to $n$. Therefore we can conclude for a general continuous spectral density by step 1).\eop

\begin{rems} \label{general_0} 1)  In Lemma \ref{mean1}, the dynamical system $(\Omega, \theta, \PP)$ can be replaced by any ergodic dynamical system.

2) If the spectral density is constant (i.e., when the $X_\k$'s are pairwise orthogonal), (\ref{orth1}) and (\ref{VarMulti1}) are a consequence of 
the law of large numbers for the number of self-intersections, that is ${V_n(\omega) \over C_0 n \ln n} \to 1$. 
The law of large numbers for $V_n(\omega, \p)$, $\p \not = \0$, is not needed. 
\end{rems}

\subsection{\bf Formulation of the quenched FCLT for a 2d random field}

\

Let $(Y_n(t)), t \in [0, 1])$ be a process on $(E, \mu)$ with values in the space  $C[0, 1)$ of real valued continuous functions on $[0,1]$ or in the space $D[0, 1)$ 
of right continuous real valued functions with left limits, endowed with the uniform norm. 

Let $(W(t), t \in [0, 1])$ be the Wiener process on $[0, 1]$. To show a functional limit theorem (FCLT) for $(Y_n(t)), t \in [0, 1])$, 
i.e., weak convergence to the Wiener process, it suffices to prove the two following properties (denoting by $\Longrightarrow$ the convergence in distribution):

{\it 1) Convergence of the finite dimensional distributions.} 
$$\forall \, 0=t_0 < t_1 <... < t_r =1, \ (Y_{n}(t_1), ..., Y_{n}(t_r)) \underset {n \to \infty} \Longrightarrow (W_{t_1}, ..., W_{t_r}),$$ 
a property which follows by the Cram\'er-Wold theorem \cite{CraWol36} from 
\begin{eqnarray}
&\sum_{j=1}^r a_j (Y_n(t_j) - Y_n(t_{j-1})) \Longrightarrow {\Cal N}(0, \sum_{j = 1}^r a_j^2 (t_j - t_{j-1})), \, \forall (a_j)_ {1 \leq j \leq r} \in \R. \label{Cramer0}
\end{eqnarray}

{\it 2) Tightness of the process.} The condition of tightness reads:
\begin{eqnarray}
&&\forall \varepsilon > 0, \ \lim_{ \delta \to 0} \limsup_n \mu(x \in E: \, \sup_{|t'-t| \leq \delta}|Y_n(\omega, x, t') - Y_n(\omega, x, t)| \geq \varepsilon) = 0. \label{tight0}
\end{eqnarray}
Now, let $(Z_n)$ be a 2-dimensional centered random walk with a finite moment of order 2 as in Subsection \ref{varRWSect1}.
Let $X = (X_\el)_{\el \in \Z^2}$ be a strictly stationary real random field, where the $X_\el$'s have zero mean and finite second moment. 
A quenched FCLT is satisfied by the r.f. $X $ if, for a.e. $\omega$, the functional central limit theorem holds for the process (cf. Notation (\ref{defsumRW}))
\begin{eqnarray}
 (Y_n(\omega, x, t))_{t \in [0, 1]} := \bigl({S_{[n t]}^{\omega, X}(x) \over \sqrt{n \log n}}\bigr)_{t \in [0, 1]}. \label{defYn0}
\end{eqnarray}
When  $(X_\el)$ is an iid random field, the model is the so-called random walk in random scenery (RWRS).
In the next section we consider first this independent case, before other non independent models in the last sections.

\section{\bf Independent random field} \label{indepSec}

\subsection{\bf Random walk in random scenery} \label{rwrc}

\

Let $X = (X_\el(x))_{\el \in \Z^2} = (T^\el f(x))_{\el \in \Z^2}$ be a 2 dimensional random field of centered i.i.d. real variables with $\E(X_0^2) = 1$ 
and mean 0 on a space $(E, \mu)$. We consider the random walk in random scenery $S_{n}^{\omega, X}(x)$ and the process defined by (\ref{defYn0}).

It was shown by E. Bolthausen \cite{Bo89} that this process satisfies an annealed FCLT, that is: 
with respect to the probability $\PP \times \mu$, the law of $Y_n$ converges weakly to the Wiener measure. 

A quenched FCLT under the assumption $\E [|X_\0|^2 \, (\log_+ |X_\0|)^\chi] < \infty$, for some $\chi > 0$, has been proved for $(Y_n(\omega, x, t) )$
in \cite{GuPoRS14}, based on \cite{Bo89}, a result of E. Bolthausen and A-S. Sznitman (2002) and a truncation argument. 

In this section, we give a direct proof of the quenched FCLT for an iid r.f. (and for moving averages of an iid r.f. in Section \ref{movAverag}), 
assuming only the existence of a moment of order 2 for the r.f.
As in \cite{Bo89} for the annealed FCLT, our proof follows the method of Newman and Wright \cite{NeWr81} for associated r.v.s.
\begin{definition} {\rm (cf. \cite{EPW})
Recall that real random variables $X_1,\ldots, X_n$ are {\it associated} if, for every $n$, for all non-decreasing (in each coordinate) functions $f,g: \R^n\mapsto \R$,
we have, if the covariance exists: ${\rm Cov}(f(X_1,\ldots,X_n), g(X_1,\ldots, X_n))\ge 0$.
A collection $X_1, X_2,\ldots$ of variables is said to be associated if every finite sub-collection is associated.

It is known that every subset of an associated family is associated. Moreover, every collection of non-decreasing functions of a family of associated random variables are associated. 
It follows that if $(X_k)$ is an associated family, in particular independent, then $(X_{Z_k(\omega)})$ is an associated family for every $\omega \in \Omega$. }
\end{definition}
\begin{thm} \label{FCLTind1} If $\E(X_\0^2) = 1$, for $\PP$-a.e. $\omega$, the process 
$\displaystyle \bigl (Y_n(\omega, x, t\bigr)_{t \in [0,1]} = \bigl({S_{\lfloor nt \rfloor}^{\omega, X} (x) \over\sqrt{ n \log n}}\bigr)_{t \in [0,1]}$
satisfies a FCLT with asymptotic variance $\sigma^2 = (\pi \sqrt{\det \Sigma})^{-1}.$
\end{thm}
\proof  1) For the convergence of the finite dimensional distributions, the proof, relying on Cram\'er-Wold's theorem and Lindeberg's CLT, 
is as in Bolthausen (\cite{Bo89}). Another proof, based on truncation and cumulants, is like the more general case of moving averages in Section \ref{movAverag}.

2) {\it Tightness of the process $(Y_n)$.} The following is shown in the proof of Theorem 3 in \cite{NeWr81}:

Let $U_1, U_2, \ldots$ be centered associated random variables with finite second order moment.
Put $S_k=\sum_{j=1}^k U_j,$ for $k\geq1$. Then, for every $\lambda>0$ and $n \geq 1$, we have
\begin{eqnarray}
\mu(\max_{1\le k \le n}|S_k| \geq \lambda \, \|S_n\|_2) \leq 2\mu \bigl(|S_n|\ge (\lambda-\sqrt{2}) \, \|S_n\|_2 \bigr). \label{NWinequal1}
\end{eqnarray}
Inequality (\ref{NWinequal1}) can be applied to $U_j =  X_{Z_j(\omega)}$ for every fixed $\omega$, as well as
to the sums $S_J=\sum_{j=b}^{b+k} X_{Z_j}$ for any interval $J = [b , b+k] \subset [0, n]$. We also note that $\E(S_J^2) = \|X_0\|_2^2 \, V(\omega, J)$.  

a) First, let us assume that $\E(X_\0^4)  < \infty$. With $K$ given by (\ref{Komega}), we have 
\begin{eqnarray}
&&\|\sum_{i \in J} X_{Z_i(\omega)}\|_{4, \mu}^4 
= 3 \E(X_{\0}^2)^2 \, \sum_{\el_1 \not = \el_2} w(\omega, J, \el_1)^2 \, w(\omega, J, \el_2)^2 + \, \E(X_\0^4) \, \sum_{\el} w(\omega, J, \el)^4 \nonumber \\
&&\ \leq 4 \, \E(X_\0^4) \, V(\omega, J)^2 \leq 4 \, \E(X_\0^4) \, (K(\theta^b \omega))^2 \, (k \ln k)^2. \label{mom4iid1}
\end{eqnarray}
Let $C_1$ be a constant $> 0$ such that $\mu\{\omega: \, K(\omega) \leq C_1\} > 0$. 
Using Lemma \ref{mean1}, for $n$ big enough and $\delta \in ]0, 1[$, there are times 
$0 = \rho_1 < \rho_2< ... < \rho_v \leq n < \rho_{v+1}$, with $v < 2/\delta$, such that $K(\theta^{\rho_i} \omega) \leq C_1$ 
and $\frac12 \delta n \leq \rho_{i+1} - \rho_i \leq \frac32 \delta n$, for $i = 1, \ldots, v$. 

Let $t_i=\frac{\rho_i}{n}$, $\lambda=\frac{\varepsilon}{\sqrt{\delta}}$, $J_i=[\rho_{i-1}, \ldots, \rho_i[$, $m_i =\frac23 (\rho_{i+1} - \rho_i) \leq \delta n$. 
There is $C$ such that, by (\ref{Komega}) and (\ref{mom4iid1}),
\begin{eqnarray}
&&\|\sum_{j=\rho_{i-1}}^{\rho_i} X_{Z_j(\omega)}\|_2 \leq C \, \|X_0\|_2 \, (n \, \delta \, \log (n \delta))^\frac12, 
\|\sum_{j=\rho_{i-1}}^{\rho_i} X_{Z_j(\omega)}\|_4 \leq C \, \|X_0\|_4 \, (n \, \delta \, \log (n \delta))^\frac12, \forall i.
\end{eqnarray}
Using (\ref{NWinequal1}), we get, with $\sigma^{(i)} = \|\sum_{j=\rho_{i-1}}^{\rho_i} X_{Z_j(\omega)}\|_2$, $\lambda_i = \varepsilon \sqrt{n\log n} / \sigma^{(i)}$, by Chebyshev's inequality (for moment of order 4):
\begin{eqnarray*}
&&\mu(\sup_{t_{i-1}\le s\le t_i} | \sum_{j=[t_{i-1}n\rfloor}^{\lfloor sn\rfloor} X_{Z_j (\omega)} |\geq \varepsilon \sqrt{n\log n}) 
=\mu(\max_{\rho_{i-1} \leq k \leq \rho_i} |\sum_{j=\rho_{i-1}}^{k} X_{Z_j(\omega)}|\geq \lambda_i \, \sigma^{(i)}) \\
&&\leq 2 \mu(|\sum_{j=\rho_{i-1}}^{\rho_i} X_{Z_j(\omega)}| \geq (\lambda_i -\sqrt{2}) \, \sigma^{(i)})
\leq 2 \mu(|\sum_{j=\rho_{i-1}}^{\rho_i} X_{Z_j(\omega)}| \geq \frac12 \lambda_i \, \sigma^{(i)}) \\
&&\leq 2 \mu(|\sum_{j=\rho_{i-1}}^{\rho_i} X_{Z_j(\omega)}| \geq \frac12 \varepsilon \sqrt{n\log n}) 
\leq 2 \, C^4 \, \|X_0\|_4^4 \, {(n \, \delta \, \log (n \delta))^2 \over \frac1{16} \varepsilon^4 (n\log n)^2} 
\leq 32 \, C^4 \, \|X_0\|_4^4 \, {\delta^2 \over \varepsilon^4}.
\end{eqnarray*}
We have used that $\lambda_i$ is big if $\delta$ is small. Observe now that (cf. \cite{Bi99})
\begin{eqnarray*}
&&\mu(\sup_{|t'-t|\le\delta} |Y_n(t)-Y_n(s)| \geq 3 \varepsilon) \leq
\sum_{i=1}^v \mu(\sup_{t_{i-1}\le s\le t_i} | \sum_{j=[t_{i-1}n\rfloor}^{\lfloor sn\rfloor} X_{Z_j (\omega)} |\geq \varepsilon \sqrt{n\log n}).
\end{eqnarray*}
Hence we get
$$\mu(\sup_{|t'-t|\le\delta} |Y_n(t)-Y_n(s)| \geq 3 \varepsilon) 
\leq 32 \, C^4 \, \|X_0\|_4^4 \, \frac{2}{\delta}\frac{\delta^2}{\varepsilon^4} = 64 \, C^4 \, \|X_0\|_4^4 \, \frac{\delta}{\varepsilon^4}.$$

b) Now we use a truncation. For $L>0$, let 
\begin{eqnarray*}
&&\hat X_\k^L:=X_\k \, {\bf 1}_{\{X_\k \leq L\}}- \E(X_\k \, {\bf 1}_{\{X_\k \leq L\}}), 
\ \tilde X_\k^L:=X_\k -\hat X_\k^L = X_\k \, {\bf 1}_{\{X_\k > L\}}- \E(X_\k \, {\bf 1}_{\{X_\k > L\}}), \\
&&\hat Y_n^L(t) = \frac{1}{\sqrt{C_0n\log n}}\sum_{j=0}^{\lfloor tn \rfloor} \hat X^L_{Z_j(\omega)} \text { and }
\tilde Y^L_n(t) := Y_n(t) - \hat Y_n^L(t)= \frac{1}{\sqrt{C_0 n \log n}} \sum_{j=0}^{\lfloor tn \rfloor} \tilde X^L_{Z_j(\omega)}.
\end{eqnarray*}
Since we have still sums of associated random variables, all what we have done above (including (\ref{NWinequal1}) holds for both sums, 
except that for the unbounded part of the truncation we only have a moment of order 2. 
We use Chebyshev's inequality (for moment of order 2) to control the unbounded truncated part:
\begin{eqnarray*}
\mu(|\sum_{j=\rho_{i-1}}^{\rho_i} \tilde X^L_{Z_j(\omega)}| \geq \frac12 \varepsilon \sqrt{n\log n}) 
\leq C^2 \, \|\tilde X^L_0\|_2^2 \, {n \, \delta \, \log (n \delta) \over \frac14 \varepsilon^2 \, n\log n} 
\leq 4 C^2 \, \|\tilde X^L_0\|_2^2 \, {\delta \over \varepsilon^2}.
\end{eqnarray*}
Hence, for $n$ and $\lambda$ big enough, the sum over $i$ is comparable for some constant $C'$ with
$$C\sum_{i=1}^v \frac{ \|\tilde X^L_0\|_2^2}{\lambda^2}\leq \frac{C'}\delta \frac{\delta \|\tilde X^L_0\|_2^2 }{\varepsilon^2} 
= C' \, \varepsilon^{-2} \, \|\tilde X^L_0\|_2^2.$$
Allying the inequality $\mu(|f+g|\ge\varepsilon) \leq \mu(|f|\geq \frac\varepsilon2)+\mu(|g| \geq \frac\varepsilon2)$
to $Y_n(t) = \hat Y_n^L(t) + \tilde Y^L_n(t)$, we obtain the bound:
\begin{eqnarray*}
&&\mu (\sup_{|t'-t| \leq \delta}|Y_n(t') - Y_n(t)| \geq 3\varepsilon) \leq 16 C_1 \,\frac{L^4 \delta}{\varepsilon^4} + 4 C' \, \frac{\|\tilde X^L_0\|_2^2}{\varepsilon^2}.
\end{eqnarray*}
We need, for fixed $\varepsilon > 0$, $\lim_{\delta \to 0^+} \limsup_n \mu (\sup_{|t'-t| \leq \delta}|Y_n(t') - Y_n(t)| \geq 3\varepsilon)  = 0$.

Let $\eta > 0$. First we take $L$ such that  $4 C' \, \frac{\|\tilde X^L_0\|_2^2}{\varepsilon^2} < \frac12 \eta$, then $\delta$ such that
$16 C_1 \,\frac{L^4 \delta}{\varepsilon^4} \leq \frac12 \eta$. \eop

\subsection{\bf A model based on the Lorentz process} \label{LorentzFCLT}

\

In this subsection we sketch briefly how to obtain an version of a FCLT where the random walk is replaced by the movement of a particle in a dispersing periodic billiard.
We refer to \cite{Pen09} and \cite{Pen14} for more details on this model.

Let be given a ``billiard table''  in the plane, union of $\Z^2$-periodically distributed obstacles with pairwise disjoint closures. 
We consider a point particle moving in the complementary $Q$ of the billiard table in $\R^2$ with unit speed and elastic reflection off the obstacles. 
By sampling the flow at the successive times of impact with the obstacles, we obtain a Poincar\'e's section of the billiard flow, the billiard transformation.

We suppose (dispersing billiard) that the obstacles are strictly convex with pairwise disjoint closures and boundaries of class $C^{r+1}$ with curvature $>0$
(Sinai's billiard or Lorentz's process). Moreover we assume a finite horizon (the time between two subsequent reflections is uniformly bounded).  

Suppose that to each obstacle is associated a real random variable with zero expectation, positive and finite variance, 
independent of the motion of the particle and that the family of these r.v.s is i.i.d.

Like in an infinite ``pinball'' with random gain, at each collision with an obstacle, the particle wins the amount given by the random variable associated with the obstacle 
which is met. Let $W_n$ be the total amount won by the particle after $n$ reflections occur. An annealed FCLT for $W_n$ has been shown by F. P\`ene (\cite{Pen09}): 
there exists $\beta_0 > 0$ such that $\displaystyle {W_{[nt]} \over \beta_0 \, n \lg n}$ converges weakly to the standard Wiener process.

To extend the result to a quenched version, we use \cite[Proposition 7]{Pen09}, in place of Inequality (\ref{majwnmIneg0}) for the r.w., and \cite[Corollary 4]{Pen14} 
(the main and most difficult step), which gives a law of large numbers for the self-intersections of the billiard transformation replacing (\ref{equivAn20}). 
Then, Remarks \ref{general_0} and the preceding method for the r.w. in random sceneries yield the quenched version of the FCLT for this model.

\section{\bf Cumulants and CLT} \label{Sectcumulant}

For the models of r.f. in Sections \ref{movAverag} and \ref{algMod1}, we need to introduce some tools. In the section, we recall the method of cumulants.

The method of cumulants recalled below can be helpful to prove the CLT in dynamical systems. 
In 1960, Leonov (\cite{Leo60a}, \cite{Leo60b}) applied it to a single algebraic endomorphism of a compact abelian group. 
In \cite{CohCo16}, \cite{CohCo17}, it was applied to multidimensional actions by algebraic endomorphisms. 

{\bf Moments and cumulants} 

In this subsection, the random variables are assumed to be uniformly bounded and centered.

Let $(X_1, ... , X_r)$ be a random vector. For $I = \{i_1, ..., i_p\} \subset \{1, ..., r\}$, let $m(I) = m(i_1, ..., i_p):= \E(X_{i_1} \cdots X_{i_p})$. 

A definition of the cumulant of $(X_1, ... , X_r)$  using the moments is
\begin{eqnarray}
C(X_1, ... , X_r) = \sum_{\Pi = \{I_1, I_2, ..., I_{p(\Pi)}\} \in \Cal P} (-1)^{p(\Pi)-1} (p(\Pi) - 1)! \, \,  m(I_1) \cdots m(I_{p(\Pi)}), \label{cumForm0}
\end{eqnarray}
where $\Pi = \{I_1, I_2, ..., I_{p(\Pi)}\}$ runs through the set $\Cal P$ of partitions of $\{1, ..., r\}$ into nonempty subsets
and $p(\Pi)$ is the number of elements of $\Pi$.

For example, if $r= 4$, the cumulant of centered r.v.s is
$$C(X_1, X_2, X_3, X_4) = \E(X_1 X_2 X_3 X_4) - [\E(X_1 X_2) \, \E(X_3 X_4)  + \E(X_1 X_3) \, \E(X_2 X_4) + \E(X_1 X_4) \, \E(X_2 X_3)].$$
Putting $s(I) := C(X_{i_1}, ... , X_{i_p})$ for $I = \{i_1, ..., i_p\}$, we have
\begin{eqnarray}
\E(X_{1} \cdots X_{r}) = \sum_{\Pi =  \{I_1, I_2, ..., I_{p(\Pi)}\} \in \Cal P} s(I_1) \cdots s(I_{p(\Pi)}). \label{cumFormu1}
\end{eqnarray}
For a single random variable $Y$, the cumulant of order $r$ is defined by $C^{(r)}(Y) := C((Y, ..., Y)_r)$, where $(Y, ..., Y)_r$ is the vector with $r$ components
equal to $Y$. If $Y$ is centered, we have 
$$C^{(2)}(Y) = \|Y\|_2^2, \ C^{(4)}(Y) = \E(Y^4) - 3 \E(Y^2)^2, \ \E(Y^4) = 3 \E(Y^2)^2 + C^{(4)}(Y).$$
If $(X_\el)_{\el\in \Z^d} = (T^{\el} f)_{\el\in \Z^d}$ is a stationary random field, we put $C_f(\el_1, ... , \el_r) =C(X_{\el_1}, ... , X_{\el_r})$.
 
Let us recall a criterium in terms of cumulants for the CLT (as well as for the convergence of the normalised moments toward those of the normal law) 
(cf. \cite[Th. 7]{Leo60b}, \cite[Th. 6.2]{CohCo17}).
\begin{thm} \label{Leonv0} If $(w_n)_{n \geq 1}$ is a summation sequence on $\Z^d$ which is $\xi$-regular (cf. Subsection \ref{varSect}), the condition
\begin{eqnarray}
\sum_{(\el_1, ..., \el_r) \, \in (\Z^d)^r} \, w_n(\el_1) ... w_n(\el_r) \, C_f(\el_1, ... , \el_r) = o \bigl((\sum_{\el \in \Z^d} \, w_n^2(\el))^{\frac{r}2}\bigr),
\ \forall r \geq 3, \label{smallCumul0}\\
\text{ implies } \ \  \bigl(\sum_{\el \in \Z^d} \, w_n^2(\el)\bigr)^{-\frac12} \, \sum_{\el \in \Z^d} w_n(\el) f(T^\el .) {\underset{n \to \infty}
\Longrightarrow } \Cal N(0, \xi(\varphi_f)). \label{cvgce1}
\end{eqnarray}
\end{thm}

The following result (cf. \cite[Lemma 6.6]{CohCo17})
gives a sufficient condition for the asymptotic nullity of the cumulants.
\begin{proposition} \label{skTozeroLem} Let $(T^\el, \el \in \Z^d)$ be a $\Z^d$-measure preserving action on a probability space $(E, \mu)$. 
If it is mixing of order $r \geq 2$, then, for any $f \in L_0^\infty(X)$, 
\begin{eqnarray}
\underset{\max_{i \not = j} \|\el_i - \el_j\| \to \infty} \lim C(T^{\el_1} f,..., T^{\el_r} f) = 0. \label{lem66}
\end{eqnarray}
\end{proposition}
Remark that (\ref{lem66}) does not give the quantitative estimate needed in (\ref{smallCumul0}). Nevertheless, it will suffice in Section \ref{algMod1} for an action 
by automorphisms of a compact abelian group which is mixing of order $r \geq 2$, in particular on a torus and $f$ is a trigonometric polynomial. 

\vskip 2mm
{\it Array of sequences and finite dimensional distributions}

Using Theorem \ref{Leonv0}, we are going to deduce from the following two conditions the asymptotic normality (after normalization) of the vectorial process 
$\displaystyle \bigl(\sum_{\el \in \Z^d} \, w_{n,1}(\el) \, T^\el f, ...,  \sum_{\el \in \Z^d} \, w_{n,s}(\el) \, T^\el f \bigr)$:
\hfill \break - asymptotic orthogonality:
\begin{eqnarray}
\int_{\T^d} (\sum_{\el \in \Z^d} w_{n,j}(\el) \, e^{2\pi i \langle \el, \t\rangle}) \, (\sum_{\el \in \Z^d} w_{n,j'}(\el) \, e^{-2\pi i \langle \el, \t\rangle}) 
\, e^{-2\pi i \langle \p, \t\rangle} \, d\t \nonumber \\
\ \ = o\bigl(\sum_{\el \in \Z^d} (w_{n, j}(\el))^2 + \sum_{\el \in \Z^d} (w_{n, j'}(\el))^2\bigr), \forall j \not = j', \forall \p \in \Z^d. \label{ortho1}
\end{eqnarray}
\hfill \break - convergence to 0 of the normalized cumulants of order $\geq 3$:
\begin{eqnarray}
&&\sum_{(\el_1, ..., \el_r) \, \in (\Z^d)^r} \, w_{n, i_1}(\el_1) ... w_{n, i_r}(\el_r) \, C(X_{\el_1}, ..., X_{\el_r}) \nonumber \\
&&\ \ \ = o(\sum_{\el \in \Z^d} \, [(w_{n, 1}(\el))^2 + ... + (w_{n, s}(\el))^2])^{r/2}, \, \forall (i_1, ..., i_r) \in \{1, ..., s\}^r, \, \forall r \geq 3, \label{negligCum3}
\end{eqnarray}

\begin{proposition} \label{CltSum1} Let $(w_{n, j})_{n \geq 1}$, $j= 1, ..., s$, be $\xi_j$-regular summation sequences, 
where the $\xi_j$'s are probability measures on $\T^ d$.  Under Conditions (\ref{ortho1}) and  (\ref{negligCum3}), the vectorial process 
$$\bigl({\sum_{\el \in \Z^d} \, w_{n, 1}(\el) \, T^\el f \over (\sum_{\el \in \Z^d} w_{n, 1}(\el)^2)^{\frac12}}, ..., 
{\sum_{\el \in \Z^d} \, w_{n, s}(\el) \, T^\el f \over (\sum_{\el \in \Z^d} w_{n, s}(\el)^2)^{\frac12}}\bigr)_{n \geq 1}$$
is asymptotically distributed as ${\cal N}(0, J_s)$, where $J_s$ is the $s$-dimensional diagonal matrix with diagonal $\xi_j(\varphi_f)$.
\end{proposition}
\proof The hypothesis (\ref{ortho1}) implies 
\begin{eqnarray}
(\sum_j a_j^2 \sum_{\el \in \Z^d} (w_{n, j}(\el))^2)^{-1} |\sum_{\el \in \Z^d} \sum_j a_j w_{n, j}(\el) \, e^{2\pi i \langle \el, \t\rangle}|^2 
\, \overset {\text{weakly}} {\underset {n \to\infty} \longrightarrow} \, (\sum_j a_j^2)^{-1} \, \sum_j a_j^2 \, \xi_j. \label{sumVar1}
\end{eqnarray}
For $s$ non zero real parameters $a_1, ..., a_s$, let $(w_n^{a_1, ..., a_s})_{n \geq 1}$ be defined by  
$$w_n^{a_1, ..., a_s}(\el) = a_1 w_{n, 1}(\el) + ... + a_s w_{n, s}(\el).$$ 

By the Cram\'er-Wold theorem, for the conclusion of the theorem, it suffices to show that the process $\sum_{\el \in \Z^d} \, w_n^{a_1, ..., a_s}(\el) \, T^\el f$ 
after normalization satisfies the CLT: 
\begin{eqnarray}
{\sum_{\el \in \Z^d} \, w_n^{a_1, ..., a_s}(\el) T^\el f \over (a_1^2 \sum_{\el \in \Z^d} (w_{n, 1}(\el))^2 + ... + a_s^2  \sum_{\el \in \Z^d} \, (w_{n, s}(\el))^2)^{\frac12}}
\Longrightarrow \Cal N(0, \sum_{j=1}^s a_j^2 \,\xi_j(\varphi_f) / \sum_{j=1}^s a_i^2). \label{TCLFd1}
\end{eqnarray}
By (\ref{negligCum3}), the sum $\displaystyle \sum_{i_1, ..., i_r \in \{1, ..., s\}^r} \, \sum_{(\el_1, ..., \el_r) \, \in (\Z^d)^r} \, w_{n, i_1}(\el_1) ... w_{n, i_r}(\el_r) 
\, C(X_{\el_1}, ..., X_{\el_r})$ satisfies (\ref{smallCumul0}) and the result follows from Theorem \ref{Leonv0}. \eop

\vskip 3mm
The following lemma will be useful in the proof of the asymptotic normality for the finite dimensional distributions.

Let $(w_n)_{n \geq 1}$ be a summation sequence on $\Z^d$ which is $\xi$-regular. For $f \in L^2(\mu)$, we put $\sigma_n(f) := \|\sum_\el w_n(\el) \, T^\el f \|_2$.
\begin{lem} \label{tclRegKer} Let $f, f_k, k = 1, 2, ..., \in L^2(\mu)$ satisfying (\ref{sumcorr1}) such that $\|\varphi_{f - f_k}\|_\infty \to 0$. Then
\begin{eqnarray*}
&&\sigma_n(f_k)^{-1} \, \sum_{\el \in \Z^d} w_n(\el) \, T^\el f_k \  {\underset{n \to \infty} \Longrightarrow} \ \Cal N(0,1), \, \forall k  \geq 1, \\
\text{ implies } &&\sigma_n(f)^{-1} \, \sum_{\el \in \Z^d} w_n(\el) \, T^\el f \ {\underset{n \to \infty} \Longrightarrow} \ \Cal N(0,1). 
\end{eqnarray*}
\end{lem}
\proof Let $(\varepsilon_k)$ be a sequence of positive numbers tending to 0, such that $\|\varphi_{f - f_k}\|_\infty \leq \varepsilon_k$. 
Let us consider the processes defined respectively by
\begin{eqnarray*}
U_n^{k} := (\sum_{\el \in \Z^d} \, w_n^2(\el))^{-\frac12} \, \sum_{\el \in \Z^d} w_n(\el) \, T^\el f_k, 
\ U_n := (\sum_{\el \in \Z^d} \, w_n^2(\el))^{-\frac12} \, \sum_{\el \in \Z^d} w_n(\el) \, T^\el f.
\end{eqnarray*}
By the $\xi$-regularity of $(w_n)$, we have:
$$(\sum_{\el \in \Z^d} \, w_n^2(\el))^{-1} \, \|\sum_{\el \in \Z^d} w_n(\el) \, T^\el f \|_2^2 
= \int_{\T^d} \, \tilde w_n \, \varphi_f \, dt \underset{n \to \infty} \to \xi(\varphi_f).$$
We can suppose $\xi(\varphi_f) > 0$, since otherwise the limiting distribution is $\delta_0$. We have $\xi(\varphi_{f - f_k}) \to 0$ (cf. Remark \ref{densRem1}). 
It follows that $\xi(\varphi_{f_k}) \not = 0$ for $k$ big enough.

The hypotheses imply ${U_n^{k} \ {\underset{n \to \infty} \Longrightarrow} \Cal N(0,\xi(\varphi_{f_k}))}$ for every $k$. Moreover, since
\begin{eqnarray*}
\lim_n \int |U_n^k - U_n|_2^2 \ d\mu &=& \lim_n \int_{\T^d} \, \tilde w_n \,\varphi_{f-f_k} \, d \t = \xi(\varphi_{f - f_k}) \leq \varepsilon_k,
\end{eqnarray*}
we have $\limsup_n \mu[|U_n^k - U_n| > \delta] \leq \delta^{-2} \limsup_n \int |U_n^k - U_n|_2^2 \ d\mu \underset{k \to \infty} \to 0$, for every $\delta > 0$.

Therefore the conclusion $U_n \ {\underset{n \to \infty} \Longrightarrow } \Cal N(0, \xi(\varphi_{f}))$ follows from \cite[Theorem 3.2]{Bi99}. \eop

\section{\bf Moving averages of iid random variables} \label{movAverag}

Let $(X_\el)_{\el \in \Z^2}$ be a r.f. of centered i.i.d. real random variables such that $\|X_\0\|_2 = 1$. Let $(a_\q)_{\q\in \Z^2}$ be an array of real numbers
such that $\sum_{\q \in \Z^2} |a_\q| < \infty$ and let  $(\Xi_\el)_{\el \in \Z^2}$ be the random field defined by $\Xi_\el(x) = \sum_{\q \in \Z^2} a_\q X_{\el - \q}(x)$. 

The correlation is 
$\displaystyle \widehat \varphi_\Xi(\el) = 
\langle \sum_{\q \in \Z^2} a_\q X_{\el - \q}, \sum_{\q' \in \Z^2} a_{\q'} X_{- \q'} \rangle = \sum_{\q \in \Z^2} a_\q \, a_{\q - \el}$.  
We have 
$$\sum_\el |\widehat \varphi_\Xi(\el)| \leq \sum_\el \sum_{\q \in \Z^2} |a_\q| \, |a_{\q - \el}| = (\sum_{\q \in \Z^2} |a_\q|)^2< +\infty.$$
The continuous spectral density of the process $(\Xi_\el)_{\el \in \Z^2}$ is 
$\displaystyle \varphi_\Xi(\t)  = |\sum a_\q e^{2\pi i \langle\q, \t\rangle}|^2$.
We assume that the asymptotic variance is $>0$, a condition equivalent to $\sum_{\q \in \Z^2} a_\q \not = 0$.

Using the method of associated r.v.s we obtain a quenched FCLT for $S_{\lfloor nt \rfloor}^{\omega, \Xi}$ (cf. Notation (\ref{defsumRW})). An annealed FCLT can be shown with a proof along the same lines.
\begin{thm} \label{FCLTindFil} The process $\displaystyle \bigl({S_{\lfloor nt \rfloor}^{\omega, \Xi}(x) \over\sqrt{ n \log n}}\bigr)_{t \geq 0}$
satisfies a quenched FCLT with asymptotic variance $\sigma^2 = |\sum_{\q \in \Z^2} a_\q |^{2} (\pi \sqrt{\det \Sigma})^{-1}.$
\end{thm}
\proof 1) {\it Convergence of the finite dimensional distributions}

a) First assume  the random variables are bounded. 
Moreover suppose first that the series reduces to a finite sum $F = \sum_{\s \in S} a_\s X_\s$, where $S$ is a finite subset of $\Z^2$. 
The case of the series, $\Xi_\0 = \sum_{\s \in \Z^2} a_\s X_\s$,  will follow by an approximation argument.

We use Proposition \ref{CltSum1}.  Conditions (\ref{ortho1}) follows from  Lemma \ref{VarrwLem1}.  Let us check (\ref{negligCum3}).

For $F$, there is $M$ such that $C(T^{\el_1} F,..., T^{\el_r} F) = 0$, if $\max_{i,j} \|\el_i - \el_j\| > M$, because if $M$ is big enough, there is a random variable
$T^{\el_{i_1}} F$ which is independent from the others in the collection $T^{\el_1} F, ..., T^{\el_r} F$ (by finiteness of $S$). 

Let $\displaystyle w_n^{u}(\omega, \el) := \sum_{j=n t_u}^{nt_{u+1}} 1_{Z_j=\el} \leq w_n(\omega, \el)$.
Since $\displaystyle \sup_{\el_1, ..., \el_r} |C(T^{\el_1} F, T^{\el_2} F, ..., T^{\el_r} F)| < \infty$, we have 
\begin{eqnarray*}
&& |\sum_{\max_{i,j} \|\el_i - \el_j\| \leq M} C(T^{\el_1} F, T^{\el_2} F, ..., T^{\el_r} F)| \, 
w_n^{i_1}(\omega, \el_1) \, w_n^{i_2}(\omega, \el_2) ... w_n^{i_r}(\omega, \el_r) \\
&&\ \ \leq \sum_{\el} \sum_{\|\j_2\|, ..., \|\j_r\| \leq M, \, \j_1 =\0} |C(T^{\el} F, T^{\el + \j_2} F, ..., T^{\el + \j_r} F)| \, \prod_{k=1}^r w_n^{i_k}(\omega, \el+\j_k)\\ 
&&\ \leq C \sum_\el \sum_{\|\j _2\|, ..., \|\j_r\| \leq M, \, \j_1 =\0} \ \prod_{k=1}^r w_n^{i_k}(\omega, \el+j_k)
\ \leq C \sum_\el \sum_{\|\j _2\|, ..., \|\j_r\| \leq M, \, \j_1 =\0} \ \prod_{k=1}^r w_n(\omega, \el+j_k).
\end{eqnarray*}
The right hand side is less than a finite sum of sums of the form $\sum_{\el \in \Z^d} \prod_{k= 1}^r w_n(\omega, \el + {\j}_k)$ with $\{\j_1, ..., \j_r\} \in \Z^2$.

By (\ref{nVarep1}), for every $\varepsilon >  0$, there is $C_\varepsilon(\omega)$ a.e. finite such that 
$\sup_\el w_n(\omega, \el) \leq C_\varepsilon(\omega) \, n^\varepsilon$. For $r \geq 3$, take $\varepsilon < {r-2 \over 2(r-1)}$. We have then
$\sum_{\el \in \Z^d} \prod_{k= 1}^r w_n(\omega, \el + {\j}_k) \leq  C_\varepsilon(\omega)^{r-1} \, n^{\varepsilon(r-1)}\, n = o(n^{r/2})$
and (\ref{negligCum3}) is satisfied. 

Using Lemma \ref{tclRegKer}, the result can be extended to any sum $\sum_{\s \in S} a_\s X_\s$, with $ \sum_{\s \in S} |a_\s| < \infty$.

b) Now if we assume only the condition $\|X_\0\|_2 < \infty$, we use a truncation argument and apply again Lemma \ref{tclRegKer}.

\vskip 3mm
2) {\it Tightness} Let $a_\q^+= \max (a_\q, 0)$, $a_\q^-= \max (-a_\q, 0)$.
Observe that the random variables $\sum_{\q \in \Z^2} a_\q^+ X_{\el - \q}(x)$, for $\el \in \Z^2$, are associated, 
as well as $\sum_{\q \in \Z^2} a_\q^- X_{\el - \q}(x)$, for $\el \in \Z^2$.

Therefore tightness can be proved separately for both processes. The proof is like the proof of tightness in Theorem \ref{FCLTind1}

\section{\bf A sufficient condition for tightness} \label{SectMoricz}

We present now a method for tightness based on the 4th-moment.  It will be used for random fields generated by algebraic automorphisms.

A nonnegative function $G_0= (G_0(b, n), b, n \geq 0)$ is said to be super-additive if 
\begin{eqnarray}
G_0(b, 0) = 0 \text{ and } G_0(b,k) + G_0(b+k, \ell) \leq G_0(b,k+\ell) , \, \forall b \geq 0, \forall k, \ell \geq 1. \label{supadd1}
\end{eqnarray} 

Let $(W_k)$ be a sequence of real or complex r.v.s on a probability space $(E, \mu)$. With the notation
$$S_{b,k} = \sum_{r=b+1}^{b+k} W_r, \, M_{b,n} = \max_{1 \leq k \leq n} |S_{b,k}|,$$
we recall a result of M\'oricz as it is used here.
\begin{thm} (M\'oricz, \cite{Mo76}) \label{thmMoricz1} Suppose that there exists $G_0$ satisfying (\ref{supadd1}) such that
\begin{eqnarray}
\E_\mu(|S_{b,n}|^4) \leq G_0(b,n)^2, \, \forall b \geq 0, \forall n \geq 1. \label{hypoG1}
\end{eqnarray} 
Then, with the constant $C_{max} = (1 - 2^{-\frac14})^{-4}$,
\begin{eqnarray}
\E_\mu(|M_{b,n}|^4) \leq C_{max} \, G_0(b,n)^2,  \, \forall b \geq 0, \forall n \geq 1. \label{toprove1}
\end{eqnarray} 
\end{thm}
Let $X = (X_\el)_{\el \in \Z^2}$ be a strictly stationary real random field on a probability space $(E, \mu)$, where the $X_\el$'s have zero mean and finite second moment,
and let
$$S_n^{\omega, X}(x) = \sum_{i=0}^{n-1} X_{Z_i(\omega)}(x), \ S_J^{\omega}(x) =\sum_{i \in J}^{n-1} X_{Z_i(\omega)}(x).$$
The maximal inequality (\ref{toprove1}) gives a criterium of tightness for the sums along a random walk:
\begin{proposition} \label{rwFCLTB} Let $G(\omega, ., .)$, $H(\omega, ., .)$ be super-additive functions such that for a parameter $\gamma$ 
and $K_1(\omega), K_2(\omega)$ a.e. finite functions on $(\Omega, \PP)$,
\begin{eqnarray}
G(\omega, b, k) \leq  K_1(\theta^b \omega) \, k \ln k, \ H(\omega, b, k) \leq  K_2(\theta^b \omega) \, k \, (\ln k)^{\gamma}, \ G(\omega, b, k) \geq k. \label{majGH1}
\end{eqnarray}  
Suppose that the r.v.s $X_\el$ are bounded and satisfy
\begin{eqnarray}
&&\E_\mu(|S_{J}^{\omega, X}|^4) \leq G(\omega, b, k)^2 +  n^\frac12 \, (\ln n)^{-(\gamma + 1)}\, H(\omega, b, k), \forall J= [b, b+k] \subset [1, n], 
\text{ for a.e. } \omega. \label{hypoG1b}
\end{eqnarray} 
Then, for every $\delta \in ]0, 1]$ there is an integer $N(\delta)$ such that, for every $\varepsilon \in ]0, 1]$, 
$Y_n(\omega, x, t) = {1 \over \sqrt {n \ln n}} \sum_{j=1}^{[nt]} \, X_{Z_j(\omega)}$ satisfies for $n \geq N(\delta)$:
\begin{eqnarray}
&&\mu(x \in E: \, \sup_{|t'-t| \leq \delta}|Y_n(\omega, x, t') - Y_n(\omega, x, t)| \geq \varepsilon) \leq \varepsilon^{-4} \, \delta, \text{ for a.e. } \omega. \label{pretight1}
\end{eqnarray}
\end{proposition}
\proof 1) Let $c \geq 0$, $\Delta_n = n^\frac12 (\ln n)^{-2}$, $\nu= \nu_n \geq \Delta_n$,
$L_n = [{\nu_n \over \Delta_n}]$, $\nu' = \nu_n' = [{\nu_n \over \Delta_n}] \Delta_n + \Delta_n - 1$.

The integer $\nu_n$ will be chosen of order $\delta n$. We can write, with the convention that $\sum_{r=0}^{-1} =0$: 
\begin{flalign*}
&\max_{0\leq k \leq \nu}|\sum_{j=0}^{k} X_{Z_{j + c}(\omega)}| \leq \max_{0\leq k \leq \nu'}|\sum_{j=0}^{k} X_{Z_{j + c}(\omega)}|
=\max_{0 \leq u \leq[{\nu\over \Delta_n}], 1 \leq k \leq \Delta_n -1}|\sum_{r=0}^{u-1} 
\sum_{j=r\Delta_n}^{(r+1)\Delta_n -1} X_{Z_{j + c}(\omega)}+ \sum_{j = u \Delta_n}^{u \Delta_n + k - 1} X_{Z_{j + c}(\omega)}|&\\
&\leq\max_{0 < u \leq L_n, \, 1 \leq k \leq \Delta_n -1}|\sum_{r=0}^{u-1} \, \sum_{j=r\Delta_n}^{(r+1)\Delta_n -1} \, X_{Z_{j + c}(\omega)}| 
+ \max_{0 \leq u \leq L_n, \, 1 \leq k \leq \Delta_n -1}\, |\sum_{j = u \Delta_n}^{u \Delta_n + k - 1} \, X_{Z_{j + c}(\omega)}|&\\
&= \max_{0 < u \leq L_n} \, |\sum_{j=0}^{u\Delta_n -1} \, X_{Z_{j + c}(\omega)}| 
+ \max_{0 \leq u \leq L_n, \, 1 \leq k \leq \Delta_n -1}\, |\sum_{j = u \Delta_n}^{u \Delta_n + k - 1} \, X_{Z_{j + c}(\omega)}| = \hat A_n + \tilde A_n.&
\end{flalign*}

With $\hat A_n$ and $\tilde A_n$ respectively the first and the second term above, this implies 
\begin{eqnarray}
&&\mu(\max_{0\leq k \leq \nu}|\sum_{j=0}^{k} \, X_{Z_{j + c}(\omega)}|  \geq \varepsilon \sqrt {n \ln n}) 
\leq \mu(\hat A_n \geq \frac12\varepsilon \sqrt {n \ln n}) + \mu(\tilde A_n \geq \frac12\varepsilon \sqrt {n \ln n}). \label{discr1}
\end{eqnarray}
For $\tilde A_n$, since the $X_\el$'s are bounded (uniformly in $\el$ by strict stationarity), (it suffices that 
$\|\max_{1 \leq k \leq \Delta_n -1}\, |\sum_{j = u \Delta_n}^{u \Delta_n + k - 1} \, X_{Z_{j + c}(\omega)}|\|_1$ is bounded uniformly in $u$,
there is $N_1(\delta)$ such that $\mu(\tilde A_n \geq \frac12\varepsilon \sqrt {n \ln n})  \leq \varepsilon^{-1} \delta$, for $n \geq N_1(\delta)$.

For $\hat A_n$ we will apply Theorem \ref{thmMoricz1} to $W_r = \sum_{j=r\Delta_n}^{(r+1)\Delta_n -1} \, X_{Z_{c+j}(\omega)}$, with
\begin{eqnarray}
G_0(b, k) := G(\omega, c+b \Delta_n, k \Delta_n) + (\ln n)^{- (\gamma - 1)} \, H(\omega, c + b \Delta_n, k \Delta_n), \label{defG0}
\end{eqnarray}
Since $G(\omega, c+b \Delta_n, k\Delta_n) \geq G(\omega, c+b \Delta_n, \Delta_n) \geq \Delta_n =  n^\frac12 (\ln n)^{-2}$, we have
\begin{eqnarray*}
&&G^2(\omega, c+b \Delta_n,  k\Delta_n) + n^\frac12 \, (\ln n)^{-(\gamma + 1)}\, H(\omega, c+b \Delta_n, k\Delta_n) 
\leq (G_0(b, k) )^2,
\end{eqnarray*} 
for every interval $[b,  b+k\Delta_n[$. Therefore
\begin{eqnarray*}
&&\E_\mu(|\sum_{r=b+1}^{b+k} \, W_r|^4)
= \E_\mu(|\sum_{j=(b+1)\Delta_n}^{(b+k+1)\Delta_n -1} \, X_{Z_{c+j}(\omega)}|^4) \leq G_0(b, k)^2, \, \forall b \geq 0, \forall n \geq 1,
\end{eqnarray*}
which implies by (\ref{toprove1}) of Theorem \ref{thmMoricz1}: 
\begin{eqnarray*}
&&\E_\mu(\max_{1 \leq k \leq p} |\sum_{j=(b+1)\Delta_n}^{(b+k+1)\Delta_n -1} \, X_{Z_{c+j}(\omega)}|^4)
\leq C_{\max} \, G_0(b, p)^2,  \, \forall b \geq 0, \forall p \geq 1.
\end{eqnarray*} 
Putting $K(\omega) := \max(K_1(\omega), K_2(\omega))$ and using (\ref{majGH1}), we get the bound
\begin{eqnarray}
&&\|\max_{u = 1}^{L_n} |\sum_{j=0}^{u\Delta_n} \, X_{Z_{c+j}(\omega)}| \|_4^4 
\leq C_{max} \, [G(\omega, c, L_n \Delta_n) + (\ln n)^{-\gamma + 1} \, H(\omega, c, L_n \Delta_n)]^2 \nonumber \\
&&\leq  C_{max} \, K(\theta^c \omega)^2 \, [L_n \Delta_n \ln (L_n \Delta_n) + (\ln n)^{-(\gamma - 1)} \, (L_n \Delta_n) (\ln (L_n \Delta_n))^{\gamma}]^2. \label{upBound1}
\end{eqnarray}

2) For $M > 0$ big enough, the set $\Omega_M := \{\omega: \, K(\omega) \leq M\}$ has a probability $\PP(\Omega_M) \geq \frac12$. 
We apply Lemma \ref{mean1} to $\Omega_M$. There is $N_2(\delta)$ such that for $n \geq N_2(\delta)$, 
we can find a sequence $0 = \rho_{1, n} < \rho_{2, n} < ... < \rho_{v, n} \leq n < \rho_{v+1, n}$ of visit times of $\theta^k \omega$ in $\Omega_M$ 
under the iteration of the shift $\theta$, such that $\frac12 \delta n \leq \rho_{i+1, n} - \rho_{i, n} \leq \frac32 \delta n$ and $v < 2/\delta$.
By construction, $K(\theta^{\rho_{i, n}}\omega) \leq M, \forall i$.

With $c= \rho_{i, n}$, $\nu_n=\nu_{i, n} = \rho_{i+1, n} - \rho_{i, n} \leq \frac32 \delta n$, $L_{i, n} = [{\nu_{i, n} \over \Delta_n}]$
(so that $L_{i, n} \Delta_n \leq \delta n$), we deduce from the upper bound (\ref{upBound1}) 
(for $n$ big enough and using $0 \leq \ln (\delta n) \leq \ln n$, if $n \geq \delta^{-1}$):
\begin{eqnarray*}
&&\|\max_{u = 1}^{L_{i, n}} |\sum_{j=0}^{u\Delta_n} \, X_{Z_{\rho_{i, n}+j}(\omega)}| \|_4^4 \leq C_{max} M \, [\nu_{i, n} \, \ln \nu_{i, n} 
+ (\ln n)^{-\gamma + 1} \, \nu_{i, n} \, (\ln \nu_{i, n})^{\gamma}]^2\\
&&\leq C_{max} M \, [\frac32\delta n \ln (\delta n) + (\ln n)^{-\gamma + 1} \, \delta n (\ln (\delta n))^{\gamma}]^2 \leq C_{max} M \, [\frac52\delta n \ln n]^2.
\end{eqnarray*}
This implies for $\hat A_n$ (cf. (\ref{discr1})), for $i= 1,  ..., v$, for a constant $C$:
\begin{flalign*}
&\mu(\max_{0 \leq u \leq L_{i, n}}|\sum_{j=0}^{u\Delta_n -1} \, X_{Z_{j + \rho_i}(\omega)}| \geq \frac12 \varepsilon \sqrt {n \ln n})
\leq {2 C_{max} M \, (\frac52 \delta n \ln n)^2 \over (\frac12 \varepsilon \sqrt {n \ln n})^4} \leq C \, \varepsilon^{-4} \, \delta^2.&
\end{flalign*}
Putting $t_i = \rho_{i, n} /n$, we obtain (\ref{pretight1}), i.e., for $n \geq N(\delta)$ with $N(\delta)$ big enough,
\begin{flalign*}
&\mu (\sup_{|t'-t| \leq \delta}|Y_n(t') - Y_n(t)| \geq 3\varepsilon) \leq \sum_{i=1}^v \mu(\sup_{t_{i-1}\leq s \le t_i}| Y_n(s) - Y_n(t_{i-1})| \geq \varepsilon)
\leq 2 C \, \varepsilon^{-4} \, \delta^2 \, v \leq 2 C \, {\delta \over \varepsilon^4}. \eop&
\end{flalign*}

\begin{rem} \label{remSum1} 1) Let be given for $s$ in a set of indices $S$ a process $X^s = (X_\el^s)_{\el \in \Z^2}$ satisfying the hypotheses of the proposition 
 for each $s$, with the same uniform bound and the same $G, H, \gamma$. Then, if $X_\el = \sum_s a_s X_\el^s$ with $\sum_s |a_s| \leq 1$, 
the r.f. $X = (X_\el)$ satisfies these conditions of the proposition and therefore the conclusion (\ref{pretight1}).

This follows from Minkowski inequality. We have for (\ref{hypoG1b}):
\begin{eqnarray*}
\|S_{J}^{\omega, X}\|_4^4 &\leq& (\sum_s |a_s| \|S_{J}^{\omega, {X^s}}\|_4)^4 \leq (\sum_s |a_s| [G(\omega, b, k)^2 
+ n^\frac12 \, (\ln n)^{-(\gamma + 1)} \ H(\omega, b, k)]^\frac14)^4\\
&=& (\sum_s |a_s|)^4 [G(\omega, b, k)^2 + n^\frac12 \, (\ln n)^{-(\gamma + 1)} \, H(\omega, b, k)].
\end{eqnarray*}
2) If the r.v.s $X_\el$ are not bounded, the conclusion of the proposition holds under the condition 
\begin{eqnarray}
&&\lim_{n \to \infty} \mu \bigl({1 \over \sqrt{ n \ln n}} \max_{u=0}^{\lfloor {n \over \Delta_n}\rfloor} \max_{k=1}^{\Delta_n} 
\, |\sum_{j = u \Delta_n}^{u \Delta_n+ k - 1} \, X_{Z_{j}(\omega)}| \geq \varepsilon\bigr) = 0,
\text{ with } \Delta_n = {n^\frac12 \over (\log n)^2}. \label{majAn100}
\end{eqnarray}

\end{rem}

\section{\bf Algebraic models \label{algMod1}} 

We consider a second type of example, generated by the action of two commuting automorphisms on tori.
For the tightness, in this example we use Proposition \ref{rwFCLTB}. 
This method could be used also in the independent model, but with a strengthening of the moment hypothesis. 

\vskip 2mm
\subsection{\bf Algebraic actions, automorphisms of the torus}

\

{\bf $\N^d$-actions by endomorphisms on a compact abelian group}

Let $G$ be a compact abelian group with Haar measure $\mu$. The group of characters of $G$ is denoted by $\hat G$ or $H$ and the set
of non trivial characters by $\hat G^*$ or $H^*$. The Fourier coefficients of a function $f$ in $L^1(G, \mu)$ are $c_f(\chi) :=
\int_G \, \overline \chi \, f \, d\mu$, $\chi \in \hat G$.

Every surjective endomorphism $B$ of $G$ defines a measure preserving transformation on $(G, \mu)$ and a dual injective
endomorphism on $\hat G$. For simplicity, we use the same notation for the actions on $G$ and on $\hat G$.

Let $(T_1, ..., T_d)$ be a finite family of $d$ commuting surjective endomorphisms of $G$ and $T^\el = T_1^{\ell_1}... T_1^{\ell_1}$, 
for $\el = (\ell_1, ..., \ell_d) \in \Z^d$. We obtain a $\Z^d$-action $(T^\el, \el \in \Z^d)$ on $G$, which is totally ergodic if and only if the dual action is free. 
The composition with a function $f$ defined on $G$ is denoted $T^\el f$. 

{\it We assume that the action on $G$ is mixing of all orders (this holds if it is totally ergodic and $G$ is connected, which is the case of a totally ergodic action on a torus).}

Let $AC_0(G)$ denote the class of real functions on $G$ with {\it absolutely convergent Fourier series} and $\mu(f) =0$, endowed with
the norm: $\|f \|_c := \sum_{\chi \in \hat G} |c_f(\chi)| < +\infty$. 
\begin{prop} \label{resum} If $f$ is in $AC_0(G)$, the spectral density $\varphi_f$ is continuous on $\T^\rho$ and $\|\varphi_{f}\|_\infty \leq \|f\|_c^2$. 
For every $\varepsilon > 0$ there is a trigonometric polynomial $P$ such that $\|\varphi_{f - P}\|_\infty \leq \varepsilon$.
\end{prop}
\proof \ Since the characters $T^\el \chi$ for $\el \in \Z^d$ are pairwise distinct, we have the inequalities
$$\sum_{\el \in \Z^d} |\langle T^{\el}f, f\rangle| \leq \sum_{\el \in \Z^d} \sum_{\chi \in \hat G} |c_f(T^\el \chi)| \, |c_f(\chi)|
\leq \sum_{\chi \in \hat G} (\sum_{\el \in \Z^d} |c_f(T^\el \chi)|) \, |c_f(\chi)| \leq (\sum_{\chi \in \hat G} |c_f(\chi)|)^2.$$
Therefore, if $f$ is in $AC_0(G)$, then $\sum_{\el \in \Z^d} |\langle T^{\el}f, f\rangle| < \infty$, the spectral density is continuous and 
$\|\varphi_f\|_\infty \leq \varepsilon$. By this inequality, we can take for $P$ the restriction of the Fourier series
of $f$ to a finite set  $\Cal E$ in $\hat G$, where $\cal E$ is such that
$\|\varphi_{f - P}\|_\infty \leq (\sum_{\chi \in \hat G \setminus {{\Cal E}}} |c_f(\chi)|)^2 \leq \varepsilon$. \eop

\vskip 3mm
{\bf Matrices and automorphisms of the torus}

Now we will restrict to the special case of matrices and endomorphisms of the torus $G = \T^\rho$.

Every $A$ in the semigroup ${\Cal M}^*(\rho, \Z)$ of non singular $\rho \times \rho$ matrices with coefficients in $\Z$ defines a
surjective endomorphism of $\T^\rho$ and a measure preserving transformation on $(\T^\rho, \mu)$. It defines also a dual
endomorphism of the group of characters $H = \widehat {\T^\rho}$ identified with $\Z^\rho$ (this is the action by the transposed of $B$, but since
we compose commuting matrices, for simplicity we do not write the transposition). When $A$ is in the group $GL(\rho, \Z)$ of matrices
with coefficients in $\Z$ and determinant $\pm 1$, it defines an automorphism of $\T^\rho$.
Recall that $A \in {\Cal M}^*(\rho, \Z)$ acts ergodically on $(\T^\rho, \mu)$ if and only if $A$ has no eigenvalue root of unity. 

For our purpose, we consider the case of automorphisms and $d= 2$. Let $(A_1, A_2)$ be two commuting matrices in $GL(\rho, \Z)$ with determinant $\pm 1$
and $A^\el = A_1^{\ell_1} \, A_2^{\ell_2}$, for $\el = (\ell_1, \ell_2) \in \Z^2$. It defines a $\Z^2$-action $(A^\el, \el \in \Z^2)$ on $(\T^\rho, \mu)$, 
which is totally ergodic if and only if $A^\el$ has no eigenvalue root of unity for $\el \not = \0$.

The composition with a function $f$ defined on $\T^\rho$ will be denoted $f \circ A^\el$, or $A^\el f$ or $T^\el f$.

Denote by $\tilde A_j$ the corresponding linear operators on $\C^\rho$. Let $\alpha_{i,j}, i=1, ..., \rho$, be the set of eigenvalues of $A_j$ (with multiplicity). 
We write $\al_i^{\el}$ for $\alpha_{i,1}^{\ell^1} \, \alpha_{i,2}^{\ell^2}$, if $\el = (\ell^1, \ell^2)$.

Explicit examples of totally such ergodic $\Z^2$-actions can be computed like the example below 
(cf. the book of H. Cohen on computational algebraic number theory \cite{Co93}):

$A_1 = \left(
\begin{matrix} -3 & -3 & 1 \cr 10 & 9 & -3 \cr -30 & -26 & 9 \cr
\end{matrix} \right), \ \ A_2  = \left(
\begin{matrix} 11 & 1 & -1 \cr -10 & -1 & 1 \cr 10 & 2 & -1 \cr
\end{matrix}
\right).$ 

\vskip 3mm
{\it Spectral density and rate of decorrelation for automorphisms of the torus}

A sufficient condition for $f$ to be in $AC_0(\T^\rho)$ is the following decay of its Fourier coefficients:
\begin{eqnarray}
|\hat f(\k)| = O(\|k\|^{-\beta}), \text{ with } \beta > \rho. \label{regFour2b2}
\end{eqnarray}
For compact abelian groups which are connected (\cite{CohCo17}) or which belong to a special family of non connected groups (\cite{CohCo16}),
a CLT has been shown for summation either over sets or along a random walk. 
Our aim is to extend this last result at least in the case of automorphisms of a torus to a functional CLT.

{\bf Number of solutions}

We use the following result on S-unit equations (\cite{Schl90}):
\begin{thm} \label{EvScSc} (\cite[Th. 1.1]{EvScSc02}) \ Let $K$ be an algebraically closed field of characteristic 0 and
For $r \geq 2$, let $\Gamma_r$ be a subgroup of the multiplicative group $(K^*)^r$ of finite rank $\rho$.
For any $(a_1,..., a_r) \in (K^*)^r$, the number $N(a_1,..., a_r, \Gamma)$ of solutions $x = (x_1,... , x_r) \in \Gamma_r$ of the equation
\begin{eqnarray*}
a_1 x_1 + ... + a_r x_r = 1,
\end{eqnarray*}
such that no proper subsum of $a_1 x_1 + ... + a_r x_r$ vanishes, is finite and satisfies the estimate
$$N(a_1,..., a_r, \Gamma) \leq \exp ((6r)^{3r} (\rho + 1)).$$
\end{thm}

There is a decomposition of $E = \C^\rho$ into vectorial subspaces $\C^\rho = \oplus_k E_k$ which are simultaneously invariant by $\tilde A_i$, $i=1, ..., d$, 
and such that  there is a basis $B_k$ in which $\tilde A_i$ restricted to $E_k$ is represented in a triangular form with an eigenvalue 
of $\tilde A_i$ on the diagonal.

This follows easily from the fact that the commuting matrices $A_i$ have a common non trivial space $W$ of eigenvectors,
and then from an induction on the dimension of the vector space, applying the induction hypothesis to the action of the maps $\tilde A_i$ on the quotient $E/W$.

Let us now consider on the torus $\T^\rho$ a character $\chi_\gamma$, $x \to \exp(2 \pi i \langle \gamma, x\rangle)$,
where $\gamma \in \Z^\rho \setminus \{0\}$. There is $k_0$ such that the component $\gamma_0$ of $\gamma$ in $E_{k_0}$ is $\not = 0$.
Let $\delta_0$ be the dimension of $E_{k_0}$. In the basis $B_{k_0} = \{e_{k_0, 1}, ..., e_{k_0, \delta_0}\}$ of $E_{k_0}$, we denote the coordinates 
of $\gamma_0$ by $(\gamma_0^1, ..., \gamma_0^{\delta_0})$. There is $\delta_0' \in \{1, ..., \delta_0\}$ such that
$\gamma_0^i = 0, \, \forall i < \delta_0', \text{ and } v_0:=\gamma_0^{\delta_0'} \not = 0$.

Due to the triangular form, for $j = 1, 2$, we have $A_j^\ell \gamma_0 = \alpha_{k_0, j}^\ell v_0 + w(j, \ell)$, $\forall \ell \in \Z$, 
where $\alpha_{k_0, j}$ is an eigenvalue of $A_j$ and where $w(j, \ell)$ belongs to the subspace generated by $\{e_{k_0, \delta_0' +1}, ..., e_{k_0, \delta_0}\}$.

By the total ergodicity of the action, we can choose $E_{k_0}$ such that the map $\el = (\ell^1, \ell^2) \to \alpha_{k_0, 1}^{\ell^1} \, \alpha_{k_0, 2}^{\ell^2}$ is injective.

We will apply Theorem \ref{EvScSc} to the multiplicative group generated by $\alpha_{k_0, j}$, $j=1, 2$.

\vskip 3mm 
\subsection{\bf Random walks and quenched CLT} \label{rwSect}

\

Our aim is to replace the model of i.i.d. variables $(X(\el), \el \in \Z^2)$ discussed in Section \ref{indepSec} by 
$X_\el = A^{\el} f = f \circ A^\el, \, \el \in \Z^2$
generated by an observable $f$ on a torus $\T^\rho$ under the action of commuting automorphisms. 

More precisely, we consider $\el\mapsto A^{\underline \ell}$ a totally ergodic $\mathbb Z^2$-action by algebraic automorphisms of $\mathbb T^\rho$, $\rho>1$, 
defined by commuting $\rho\times \rho$ matrices $A_1, A_2$ with integer entries, determinant $\pm 1$ 
such that the eigenvalues of $A^{\el} = A_1^{\ell_1} A_2^{\ell_2}$ are $\not = 1$, if $\el = (\ell_1, \ell_2) \not = (0, 0$).

Recall the notation, for $f$ a real function on $G= \T^\rho$, $S_{n}^n(f) :=\sum_{k=1}^{n}A^{Z_k(\omega)} f$.

The following quenched FCLT extends the CLT proved in (\cite{CohCo17}).
\begin{thm} \label{rwFCLTB2} If $(Z_n)$ is a 2-dimensional reduced centered random walk with a finite moment of order 2 and $f$ is in $AC_0(\T^\rho)$ with spectral density 
$\varphi_f$ and a non zero asymptotic variance, then for a.e. $\omega$ the process
$\Big(\frac1{\sqrt{n\log n}}S_{\lfloor nt\rfloor}^\omega(f)\Big)_{t\in[0,1]}$ satisfies a FCLT holds.
\end{thm}
\proof 
{\it 1) Convergence of the finite dimensional distributions}

1a) First suppose $f$ is a trigonometric polynomial.  Let $f = \sum_{\k \in \Lambda} c_\k(f) \, \chi_\k$, where $(\chi_\k, \k \in \Lambda)$ 
is a finite set of characters on $\T^\rho$, $\chi_{0}$ the trivial character. 

We use Proposition \ref{CltSum1}: (\ref{ortho1}) follows from (\ref{orth1}) and Lemma \ref{VarrwLem1}; for (\ref{negligCum3}), we have to show that, for a.e. $\omega$,
\begin{eqnarray}
&&C_r(\sum_{\el \in \Z^2} \, w_n(\omega, \el) \, T^\el f) = o((n \ln n)^{r/2}), \, \forall r \geq 3. \label{cumulBdr0}
\end{eqnarray}
We apply Theorem \ref{Leonv0}. Let us check (\ref{smallCumul0}). For $r$ fixed, the function 
$(\n_1, ..., \n_r) \to m_f(\n_1, ..., \n_r) := \int_X \, T^{\n_1} f \cdots T^{\n_r} f \, d\mu$ takes a finite number of values, 
since $m_f$ is a sum with coefficients 0 or 1 of the products $c_{k_1} ... c_{k_r}$ with $k_j$ in a finite set. 
The cumulants of a given order take also a finite number of values according to (\ref{cumFormu1}).

Therefore, since mixing of all orders implies $\underset{\max_{i,j} \|\el_i - \el_j\| \to \infty} \lim \, C(T^{\el_1} f,..., T^{\el_r} f)  = 0$
by Proposition \ref{skTozeroLem}, there is $M_r$ such that $C(T^{\el_1} f,..., T^{\el_r} f) = 0$, if $\max_{i,j} \|\el_i - \el_j\| > M_r$.

1b) For $f \in AC_0(\T^\rho)$, using Proposition \ref{resum} and Lemma \ref{tclRegKer}, the convergence follows by approximation 
of $f$ by trigonometric polynomials $f_L$ in such a way that $\lim_L \varphi_{f - f_L} (0)= 0$.

{\it 2) Moment of order 4 and tightness}

Let us consider real centered functions $f$ on $\T^\rho$ in $AC_0(\T^\rho)$, i.e., in real form such that 
$$f(x) = \sum_{\v \not = \0} \, [a_\v(f) \cos 2 \pi \langle \v, x\rangle + b_\v(f) \sin 2 \pi \langle \v, x\rangle], \text{ with } \sum_{\v} \, [|a_v(f)| + |b(\v)|] < +\infty.$$ 
Taking into account Remark \ref{remSum1}, it suffices to consider for $f$ a character and show that the bound is uniform, independent of the character.

So we consider on the torus $\T^\rho$ a character $\chi_v: x \to \exp(2 \pi i \langle v, x\rangle)$, where $v \in \Z^\rho \setminus \{0\}$.

For an interval $J =[b, b+k] \subset [1, n]$, we have: 
\begin{eqnarray*}
&\|\sum_{i= b }^{b+k} A^{Z_i} \chi_v \|_4^4 = 
\#\{(i_1, i_2, i_3, i_4) \in J^4: \ (A^{Z_{i_1}} -A^{Z_{i_2}}+A^{Z_{i_3}}-A^{Z_{i_4}}) \, v =0 \}.
\end{eqnarray*}
This number is bounded by
$\#\{(i_1, i_2, i_3, i_4) \in J^4: \ (\al_u^{Z_{i_1}} - \al_u^{Z_{i_2}} + \al_u^{Z_{i_3}}  - \al_u^{Z_{i_4}}) v_0 = 0\}$,
where $\v_0$ is some non zero component of $\v$ in a suitable basis in which $A_1, A_2$ have a simultaneous triangular representation  (cf. the previous subsection).
The notation is: $\al_u = \alpha_{u, 1} \, \alpha_{u, 2}$, with $\alpha_{u, 1}$ (resp. $\alpha_{u, 2}$) an eigenvalue of $A_1$ (resp. $A_2$) and
$\al_u^{Z_i} = \alpha_{u, 1}^{Z_i^1} \, \alpha_{u, 2}^{Z_i^2}$.

This number is less than $G^2(\omega, b , k) + H(\omega, b , k)$, where
\begin{eqnarray*}
G(\omega, b , k) &&:= \#\{(i_1, i_2) \in J^2: \ \al_u^{Z_{i_1}} - \al_u^{Z_{i_2}} = 0 \},\\
H(\omega, b , k) &&:= \#\{(i_1, i_2, i_3, i_4) \in J^4: \, \al_u^{Z_{i_1}} - \al_u^{Z_{i_2}} + \al_u^{Z_{i_3}} - \al_u^{Z_{i_4}}= 0 \},
\end{eqnarray*}
where above in $H$ we count the number of solutions without vanishing proper sub-sums.

By the choice of the component $\v_0$, if $\al_u^{Z_{i_1}}v(\chi) = \al_u^{Z_{i_2}}v(\chi)$, then $Z_{i_1} = Z_{i_2}$. 
Therefore, $G$ is the number of self-intersections of the r.w. starting from $b$:
$$G(\omega, b , k) = \#\{(i_1, i_2) \in J^2 :\, {Z_{i_1}}={Z_{i_2}} \}.$$
For $H$, up to a permutation of indices, we can assume that $i_4 < i_3 < i_2 < i_1$. We may write up to a constant factor:
\begin{eqnarray*}
&&H(\omega, b , k)  = \#\{b \leq i_4 < i_3 < i_2 < i_1 \leq b+k : \ \al_u^{Z_{i_1}-Z_{i_4}} - \al_u^{Z_{i_2}-Z_{i_4}} 
+ \al_u^{Z_{i_3}-Z_{i_4}} =  1\}.
\end{eqnarray*}
By Theorem \ref{EvScSc}, the set of triples $\el_1, \el_2, \el_3 \in \Z^2$ (without vanishing proper sub-sum) solving the equation 
$\al_u^{\el_1} - \al_u^{\el_2} + \al_u^{\el_3} = 1$ is a finite set $F$. 

We can now apply (\ref{majnumb3}) in Lemma \ref{majwnm}: there exists a positive integrable function $C_3$ such that for $W_n(\omega, \el_1, \el_2, \el_3)$ 
defined by (\ref{wn123}), $W_n(\omega, \el_1, \el_2, \el_3) \leq C_3(\omega) \, n \, (\ln n)^5, \, \forall n \geq 1$.

Therefore $H(\omega, b , k) \leq (\Card \, F) \ C_3(\theta^b \omega) \, k \, (\ln k)^5$.

Remark that the bounds do not depend on the character, but only on $A_1, A_2$.

The tightness property follows now from Proposition \ref{rwFCLTB} with $\gamma = 5$. \eop

\end{document}